\documentclass{amsart}
\usepackage{amssymb}
\usepackage{amsfonts}

\newtheorem{theorem}{Theorem}
\theoremstyle{plain}

\newtheorem{corollary}{Corollary}

\newtheorem{lemma}{Lemma}

\newtheorem{proposition}{Proposition}

\numberwithin{equation}{section}
\input{tcilatex}

\begin{document}
\title[Short Title]{Virtual Endomorphisms of Nilpotent Groups}
\author{Adilson Berlatto}
\address{Departamento de Matem\'{a}tica, Universidade Federal de Mato
Grosso, ICLMA, Pontal do Araguaia, MT, Brazil }
\email{berlatto@cpd.ufmt.br}
\author{Said Sidki}
\address{Departamento de Matem\'{a}tica, Universidade de Bras\'{\i}lia, Bras%
\'{\i}lia DF 70910-600, Brazil}
\email{sidki@mat.unb.br.}
\thanks{The second author thanks Laurent Bartholdi for hospitality at \'{E}%
cole Polytechnique F\'{e}d\'{e}rale -Lausanne during November 2005 and
acknowledges support from the Brazilian CNPq and Finatec}
\date{}

\begin{abstract}
A virtual endomorphism of a group $G$ is a homomorphism $f:H\rightarrow G$
where $H$ is a subgroup of $G$ of finite index $m$. The triple $\left(
G,H,f\right) $ produces a state-closed (or, self-similar) representation $%
\varphi $ of $G$ on the $1$-rooted $m$-ary tree. This paper is a study of
properties of the image $G^{\varphi }$ when \ $G$ is nilpotent. In
particular, it is shown that if $G$ is finitely generated, torsion-free and
nilpotent then $G^{\varphi }$ has solvability degree bounded above by the
number of prime divisors of $m$.
\end{abstract}

\maketitle

\section{Introduction}

A virtual endomorphism of a group $G$ is a homomorphism $f:H\rightarrow G$
where $H$ is a subgroup of $G$ of finite index $m$. A recursive construction
using $f$ produces a so called \textit{state-closed} (or, \textit{%
self-similar}) representation of $G$ on a $1$-rooted regular $m$-ary tree.
The kernel of this representation is the maximal subgroup $K$ of $H$ which
is both normal in $G$ and is $f$-invariant, in the sense that $K^{f}\leq K$;
it is called the $f$-$core\left( H\right) $.

The notion of virtual endomorphisms of groups is not recent. It already
appeared in 1969, in M. Shub's \cite{Shub}, in connection with endomorphisms
of compact differentiable manifolds. State-closed groups were introduced in 
\cite{Sid}, justified by the fact that the Grigorchuk $2$-group, the
Gupta-Sidki $p$-groups, the affine group $\mathbb{Z}^{n}GL(n,\mathbb{Z)}$ 
\cite{BruSid}, as well as an automata group of Aleshin- claimed to be free
in \cite{Alesh}- satisfied such a condition. State-closed representations of
groups on the binary tree were studied in some depth in \cite{NekSid} and
dynamical aspects of these were developed by Nekrashevych into a
far-reaching theory in \cite{Nek}.

The question of existence of finite-state, state-closed representations of
certain groups, especially of free groups, stimulated a number of
interesting constructions. Glasner and Mozes \cite{GlasMoz} used ideas from
homogeneous tree lattices to obtain such a representation of a free group of
rank $14$ acting on on the $6$-tree. This was followed by a construction by
Muntyan and Savchuk (see \cite{Nek} 1.10.3) for the free group of rank $2$
on the $6$-tree. More recently, Vorobets and Vorobets \cite{Vorobets} have
produced a proof that a group defined on the binary tree, related to the one
proposed initially by Aleshin, is indeed free of rank $3$.

The present paper extends the results on free abelian state-closed groups in 
\cite{NekSid} to finitely generated nilpotent state-closed groups. The main
emphasis though is on the subclass of torsion-free groups- following P.
Hall's notation, these are $\mathfrak{T}$-groups, or $\mathfrak{T}_{c}$%
-groups when the nilpotency class is $c$.

We refer to the data $\{G$, $H\leq G$, $f:H\rightarrow G$, $\left[ G:H\right]
=m\}$ as a \textit{triple} $\left( G,H,f\right) $ of degree $m$. If the $f$-$%
core\left( H\right) $ is the trivial subgroup then $f$ and the triple $%
(G,H,f)$ are called\textit{\ simple} and when $f$ is a simple epimorphism it
is called \textit{recurrent}. If the only $f$-invariant subgroup of $G$ is
the trivial subgroup then $f$ and the triple $(G,H,f)$ are called \textit{%
strongly simple}. To give an example of a strongly simple triple, we let $G$
be the free nilpotent group $F\left( c,d\right) $ of class $c$, freely
generated by $x_{i}$ $\left( 1\leq i\leq d\right) $, $H=\left\langle
x_{i}^{n}\text{ }\left( 1\leq i\leq d\right) \right\rangle $ where $n$ is a
fixed integer greater than $1$ and let $f$ be the extension of the map $%
x_{i}^{n}\rightarrow x_{i}$ $\left( 1\leq i\leq d\right) $. For another
example, consider the group $G$ of lower triangular matrices $\left( 
\begin{array}{ccc}
1 & 0 & 0 \\ 
a & 1 & 0 \\ 
c & b & 1%
\end{array}%
\right) $ with integer entries, $H$ its subgroup of index $4$, formed by the
matrices $\left( 
\begin{array}{ccc}
1 & 0 & 0 \\ 
2a & 1 & 0 \\ 
2c & b & 1%
\end{array}%
\right) $ and define 
\begin{equation*}
f:\left( 
\begin{array}{ccc}
1 & 0 & 0 \\ 
2a & 1 & 0 \\ 
2c & b & 1%
\end{array}%
\right) \rightarrow \left( 
\begin{array}{ccc}
1 & 0 & 0 \\ 
b & 1 & 0 \\ 
-c & a & 1%
\end{array}%
\right) \text{.}
\end{equation*}

We review state-closed groups and representations in Section 2 and
illustrate how to produce concretely state-closed nilpotent groups.
Moreover, we prove that state-closed abelian groups are small, in the sense
that their centralizers in the group of automorphisms of the tree coincide
with their topological closure in this last group.

If $G$ is an abelian group then naturally, $\ker (f)\leq f$-$core(H)$ for
any triple $\left( G,H,f\right) $. The relationship between $\ker (f)$ and
the $f$-$core(H)$ for general nilpotent groups is established in Section 4$.$

\begin{theorem}
\textit{Let }$G$\textit{\ be a nilpotent group, }$H$\textit{\ a subgroup of
finite index }$m$ in $G$\textit{, }$f\in Hom(H,G)$\textit{\ and }$L=$\textit{%
\ }$f$\textit{-}$core(H)$\textit{. Then, }$\ker (f)\leq \sqrt[H]{L}$\textit{%
, the isolator of }$L$ in $H$.\textit{\ }
\end{theorem}

A group $G$ is said to be to be \textit{compressible} provided every
subgroup $H$ of finite index in $G$ contains a subgroup $K$ isomorphic to $G$%
. It was shown by G. S. Smith in \cite{Smith} that $\mathfrak{T}_{c}$-groups
are compressible when $c\leq 2$. We extend this result to strongly simple
triples in Subsection 5.1, as follows

\begin{theorem}
Let $G$ be a $\mathfrak{T}_{c}$-group with $c\leq 2$ and let $H$ be a
subgroup of finite index in $G$. Then there exists a subgroup $K$ of finite
index in $H$, which admits a strongly simple epimorphism $f:K\rightarrow G$.
\end{theorem}

More on compressibility and co-Hopfianity questions\textit{\ }concerning $%
\mathfrak{T}_{c}$-groups can be found in \cite{Beleg}.

Given an integer $m>1$, let $l\left( m\right) $ be the number of prime
divisors of $m$ (counting multiplicities) and $a\left( m\right) $ the
largest exponent of the prime divisors of $m$. Denote by $c(G)$ the
nilpotency class of $G$, by $s\left( G\right) $ the derived length of $G$
and by $d(G)$ the minimum number of generators of $G$.

If $G$ is a finitely generated nilpotent group and $H$ a subgroup of $G$ of
index $m$ then it is well-known that $\overline{G}=\frac{G}{core\left(
H\right) }$ is finite and $c\left( \overline{G}\right) \leq a\left( m\right) 
$. We find in Subsection 5.2 such limitations for $s\left( G\right) $ and $%
c\left( G\right) $ with respect to $l\left( m\right) $.

\begin{theorem}
Let $G$ be a $\mathfrak{T}$-group and $H$ a subgroup of finite index $m$ in $%
G$. If $f:H\rightarrow G$ is simple then $s\left( G\right) \leq l\left(
m\right) $. If $f$ is strongly simple then $c\left( G\right) \leq l\left(
m\right) $.
\end{theorem}

The restriction on $G$ is striking, given that $\mathfrak{T}_{c}$-groups
have faithful finite-state representations on the binary tree ( that is, $%
m=2 $) for any $c\geq 0$ (see, \cite{NekSid}). To show that in the first
part of the theorem $s\left( G\right) $ cannot be replaced by $c\left(
G\right) $, we construct in Subsection 5.4, an ascending sequence of simple
triples $\left( G_{n},H_{n},f_{n}\right) $ where the groups $G_{n}$ are
metabelian $\mathfrak{T}$-groups with $d\left( G_{n}\right) =2$, $c\left(
G_{n}\right) =n $, $\left[ G_{n}:H_{n}\right] =4$. Using another sequence of
examples, we show that the limit in the second part of the theorem to be
satisfactory.

It is important to observe that no such limitations exist for groups of
prime power order. For let $p$ be a fixed prime number, $G$ be the $s$%
-iterated wreath product $W_{s}=\left( ((C_{p}wr..)wr)C_{p}\right) wrC_{p}$, 
$H$ its base subgroup and $\pi _{1}$ the projection of $H$ on its $1$st
coordinate. Then $\left[ G:H\right] =p$ and $\left( G,H,\pi _{1}\right) $ is
strongly simple, yet $G$ has nilpotency class $p^{s}$ and derived length $%
s+1 $.

In Section 6, we prove the following divisibility relation between indices
of subgroups

\begin{theorem}
Let $G$ be finitely generated nilpotent group, $H$ a subgroup of $G$ of
finite index $\left[ G:H\right] =m$, $f:H\rightarrow G$ a monomorphism and $%
\left[ G:H^{f}\right] =m^{\prime }$. Furthermore, let $U$ be a subgroup of $%
H $ and $V=\left\langle U,U^{f}\right\rangle $. Suppose $\left[ V:U\right]
=l,\left[ V:U^{f}\right] =l^{\prime }$ are finite. Then there exist integers 
$m_{1}|m,m_{1}^{\prime }|m^{\prime }$ such that $lm_{1}^{\prime }=l^{\prime
}m_{1}$.
\end{theorem}

As an application, we obtain

\begin{theorem}
Let $G$ be a $\mathfrak{T}$-group, $H$ a subgroup of $G$ of finite index $m$
which is a square-free integer and $f:H\rightarrow G$ a simple epimorphism.
Then $G$ is abelian.
\end{theorem}

\bigskip

The combination of conditions $a\left( m\right) =1$ and $f$ being a simple
epimorphism (recurrent) in the above theorem produced $s\left( G\right) =1$;
that is, $s\left( G\right) =a\left( m\right) $. This raises the question
about possible improvements of the bound $l(m)$ in Theorem 3, under
different types of conditions. Another question concerns the impact of the
combination state-closed and finite-state would have on $\mathfrak{T}$%
-groups. For simple triples $\left( G,H,f\right) $ where $G$ is free abelian
group of finite rank and $\left[ G:H\right] =m=2$, it was shown in \cite%
{NekSid} (see also, \cite{Nek} Sec. 2.12) that the roots of the
characteristic polynomial of $f$ lie in the interior of the unit circle.

We thank the referee for a very careful reading of our paper, for a number
of positive suggestions and for providing the explicit formula in Section
5.4.

\section{State-closed groups and representations}

Let $Y$ be a non-empty set, $P\left( Y\right) $ the group of permutations of 
$Y$ and let $\mathcal{T(}Y)$ be $1$-rooted tree indexed by the free monoid $%
Y^{\ast }$ generated by $Y$. Then the group of automorphisms $\mathcal{A}%
=Aut(\mathcal{T(}Y))$ of the tree is isomorphic to the semidirect product of 
$\mathcal{A}^{Y}$ by $P\left( Y\right) $ and under this identification, we
have the decomposition $\mathcal{A=}\left( \mathcal{A}^{Y}\right) P\left(
Y\right) $. Thus, an $\alpha \in \mathcal{A}$ is represented as $\alpha
=\left( \alpha _{y}|y\in Y\right) \alpha ^{\sigma }$ where $\alpha _{y}\in 
\mathcal{A}$ and $\sigma \in Hom\left( \mathcal{A},P\left( Y\right) \right) $%
.

The set of \textit{states} of $\alpha $ is $Q\left( \alpha \right) =\left\{
\alpha _{u}|u\in Y^{\ast }\right\} $. For $G\leq \mathcal{A}$, let $%
stab_{G}\left( i\right) $ denote the subgroup of $G$ formed by elements
which leave the $i$th level vertices fixed. The group $\mathcal{A}$ is the
inverse limit of its quotients $\frac{\mathcal{A}}{stab_{\mathcal{A}}\left(
i\right) }$ and as such becomes a topological group. Also, for $u\in Y^{\ast
}$ of length $i$, let $\pi _{u}:stab_{G}\left( i\right) \rightarrow \mathcal{%
A}$ be the projection map on the $u$th coordinate. A subgroup $G$ of $%
\mathcal{A}$ is \textit{state-closed} provided $Q\left( \alpha \right) $ is
a subset of $G$ for all $\alpha \in G$ and is\textit{\ finite-state} if $%
Q\left( \alpha \right) $ is finite for all $\alpha \in G$.

On fixing $y\in Y$, we obtain the triple $\left( G,stab_{G}\left( y\right)
,\pi _{y}\right) $. In the other direction, given a triple $\left(
G,H,f\right) $, we represent the group $G$ on the right cosets of $H$,
keeping track of the factor sets, as in Schreier's theorem. Then, we use $f$
to repeat this process down a chain of subgroups. In the limit, this
produces a state-closed representation of $G$ on a rooted tree of degree $%
\left[ G:H\right] $, as follows.

\begin{theorem}
Let $G$ be a group, $H$ a subgroup of $G$, and $Y$ a right transversal of $H$
in $G$. Let $\sigma $ be the permutation representation of $G$ on $Y$; for $%
g\in G,y\in Y$, write $g^{\sigma }:Hy\rightarrow Hyg=Hy^{\prime }$, $%
y^{\prime }\in Y$ and $y^{\prime }=\left( y\right) ^{g^{\sigma }}$. Also,
let $\mathcal{T(}Y)$ be $1$-rooted tree indexed by the free monoid $Y^{\ast
} $and let $f\in Hom\left( H,G\right) $. Then, the quadruple $\left(
G,H,Y,f\right) $ provides a representation $\varphi $ of $G$ into the
automorphism group of the tree $\mathcal{T(}Y)$, defined by%
\begin{equation*}
g^{\varphi }=\left\{ \left( yg.\left( y^{g^{\sigma }}\right) ^{-1}\right)
^{f\varphi }|y\in Y\right\} g^{\sigma }
\end{equation*}%
Furthermore, $\ker \left( \varphi \right) =f$-core$\left( H\right) $.
\end{theorem}

To dispel difficulties with the notation in the above formula, we offer a
simple example. Let $G$ be the additive group of integers $\mathbb{Z}$, $H=2%
\mathbb{Z}$, and $Y=\left\{ 0,1\right\} $; then, $\sigma :0\leftrightarrow 1$%
. Define $f:2\mathbb{Z\rightarrow Z}$ by $2n\rightarrow n$. Then, $%
1^{\varphi }=\left( 0^{\varphi },1^{\varphi }\right) \sigma $ which is none
other the binary adding machine.

The proof of the theorem is a direct extension of that of Theorem 3.1 in 
\cite{NekSid} and can be found in Section 2.5 of \cite{Nek}.

\subsection{Producing nilpotent state-closed groups}

We illustrate in this subsection how certain initial conditions about a
simple triple $\left( G,H,f\right) $ lead to an understanding of the
properties of its state-closed representations, which in turn can be used to
construct examples of such triples. Consider the following configuration:%
\begin{equation*}
\text{ }\left[ G:H\right] =4,
\end{equation*}%
\begin{eqnarray*}
z &\in &Z\left( G\right) ,\left[ \left\langle z\right\rangle :\left\langle
z\right\rangle \cap H\right] =2, \\
f &:&H\rightarrow G\text{ such that }f:z^{2}\rightarrow z\text{.}
\end{eqnarray*}%
Let $Y=\left\{ y_{1},y_{2},y_{3},y_{4}\right\} $ be a right transversal of $%
H $ in $G$. Then, we may choose $y_{1}=e,y_{2}=z$ and $y_{4}=y_{3}z$. We
identify $y_{i}$ with its subscript $i$. Thus, $z^{\sigma }$ is the
permutation $\left( 1,2\right) \left( 3,4\right) $ and $z^{\varphi }=\left(
e,z^{\varphi },e,z^{\varphi }\right) .\left( 1,2\right) \left( 3,4\right) $;
we suppress $\varphi $ from the notation and thus obtain 
\begin{equation*}
z=\left( e,z,e,z\right) .\left( 1,2\right) \left( 3,4\right) \text{.}
\end{equation*}

Let $\mathcal{A}=Aut\left( \mathcal{T(}Y)\right) ,C=C_{\mathcal{A}}\left(
z\right) $. We will show that

(i) $C$ is state-closed;

(ii) the commutator equation $\left[ \alpha ,\beta \right] =z$ has the
following particular solution in $C$: 
\begin{equation*}
\alpha =\left( \alpha ,\alpha z,\alpha ,\alpha \right) \left( 1,2\right)
,\beta =\left( z,z,z^{-1}\beta ,z^{-1}\beta \right) \left( 1,3\right) \left(
2,4\right) \text{;}
\end{equation*}%
Furthermore, The group $R=\left\langle \alpha ,\beta \right\rangle $ is
isomorphic to $F\left( 2,2\right) $, is recurrent and is finite-state ;

(iii) there exists $\kappa \in N_{C}\left( R\right) $ defined by%
\begin{equation*}
\kappa =\left( \alpha ^{3}\kappa ^{2},\alpha ^{3}\kappa ^{2},\alpha \kappa
^{2},\alpha \kappa ^{2}\right)
\end{equation*}%
such that the group $S=\left\langle \alpha ,\beta ,\kappa \right\rangle $ is
a $\mathfrak{T}_{3}$ group and has the presentation%
\begin{equation*}
\left\{ \alpha ,\beta ,\kappa |\left[ \alpha ,\beta ,\alpha \right] =\left[
\alpha ,\beta ,\beta \right] =\left[ \alpha ,\beta ,\kappa \right] =\left[
\alpha ,\kappa \right] =e,\left[ \beta ,\kappa \right] =\alpha ^{-2}\right\} 
\text{.}
\end{equation*}

Proof. It is straightforward to check that the elements of the centralizer $%
C=C_{A}\left( z\right) $ in $\mathcal{A}$ are of eight types:%
\begin{eqnarray*}
x_{1} &=&\left( h_{1},h_{1},h_{2},h_{2}\right) ,x_{2}=\left(
h_{1},h_{1}z,h_{2},h_{2}\right) \left( 1,2\right) , \\
x_{3} &=&\left( h_{1},h_{1},h_{2},h_{2}z\right) \left( 3,4\right)
,x_{4}=\left( h_{1},h_{1}z,h_{2},h_{2}z\right) \left( 1,2\right) \left(
3,4\right) , \\
x_{5} &=&\left( h_{1},h_{1},h_{2},h_{2}\right) \left( 1,3\right) \left(
2,4\right) ,x_{6}=\left( h_{1},h_{1}z,h_{2},h_{2}z\right) \left( 1,4\right)
\left( 2,3\right) , \\
x_{7} &=&\left( h_{1},h_{1},h_{2},h_{2}z\right) \left( 1,3,2,4\right)
,x_{8}=\left( h_{1},h_{1}z,h_{2},h_{2}\right) \left( 1,4,2,3\right) \text{ }
\end{eqnarray*}%
where in each case, $h_{1},h_{2}\in C$. Therefore, $C$ is also state-closed.

If $x,x^{\prime }\in C$ are such that $\left[ x,x^{\prime }\right] =z$ then
there exist $x_{2},x_{5}\in \left\langle x,x^{\prime }\right\rangle $ such
that $\left[ x_{2},x_{5}\right] =z$. Let 
\begin{equation*}
x_{2}=\left( h_{1},h_{1}z,h_{2},h_{2}\right) \left( 1,2\right) ,x_{5}=\left(
k_{1},k_{1},k_{2},k_{2}\right) \left( 1,3\right) \left( 2,4\right) .
\end{equation*}%
Then $\left[ x_{2},x_{5}\right] =z$ if and only if

\begin{equation*}
h_{2}=k_{1}^{-1}h_{1}k_{1},\left[ h_{1},k_{1}k_{2}\right] =z\text{.}
\end{equation*}%
One solution is $h_{1}=h_{2}=x_{2}$, $k_{1}=z$, $k_{2}=z^{-1}x_{5}$. With
these choices, rename $x_{2}$ as $\alpha $ and $x_{5}$ as $\beta $; thus, 
\begin{eqnarray*}
\alpha &=&\left( \alpha ,\alpha z,\alpha ,\alpha \right) \left( 1,2\right) ,
\\
\beta &=&\left( z,z,z^{-1}\beta ,z^{-1}\beta \right) \left( 1,3\right)
\left( 2,4\right) \text{.}
\end{eqnarray*}%
It can be verified directly that $\alpha ,\beta $ generate a group $R$
isomorphic to $F(2,2)$, that $R$ is finite-state, and recurrent.

We search in the normalizer of $R$ in $C$ for an element $\kappa $ such that 
$\left\langle \alpha ,\beta ,\kappa \right\rangle $ is nilpotent of class $3$%
. Clearly, we may assume $\kappa =\left( \kappa _{1},\kappa _{1},\kappa
_{2},\kappa _{2}\right) $. Since $\kappa ^{-1}\kappa ^{\alpha },\kappa
^{-1}\kappa ^{\beta }$ stabilize the $1$st level of the tree, there exist
integers $i,j,k,l,s,t$ such that 
\begin{equation*}
\kappa ^{-1}\kappa ^{\alpha }=\alpha ^{2i}\beta ^{2j}z^{r},\kappa
^{-1}\kappa ^{\beta }=\alpha ^{2l}\beta ^{2k}z^{s}
\end{equation*}%
and thus, modulo $R^{\prime }$, we have 
\begin{equation*}
\alpha ^{\kappa }=\alpha ^{1-2i}\beta ^{-2j},\beta ^{\kappa }=\alpha
^{-2l}\beta ^{1-2k}\text{.}
\end{equation*}%
Since the action of $\kappa $ on $R/R^{\prime }$ is nilpotent, the matrix $%
\left( 
\begin{array}{cc}
1-2i & -2j \\ 
-2l & 1-2k%
\end{array}%
\right) $ has determinant $1$ and trace $2$. Thus, 
\begin{equation*}
k=-i,jl=-i^{2}\text{.}
\end{equation*}

Rather than describing all possible solutions we try 
\begin{eqnarray*}
i &=&j=k=0,l=1\text{,} \\
r &=&s=0\text{.}
\end{eqnarray*}%
Then, 
\begin{eqnarray*}
\kappa &=&\left( \kappa _{1},\kappa _{1},\kappa _{2},\kappa _{2}\right) , \\
\kappa ^{\alpha } &=&\kappa ,\kappa ^{-1}\kappa ^{\beta }=\alpha ^{2}\text{.}
\end{eqnarray*}%
Now, we calculate%
\begin{eqnarray*}
\kappa ^{-1}\kappa ^{\alpha } &=&\left( \kappa _{1}^{-1}\kappa _{1}^{\alpha
},\kappa _{1}^{-1}\kappa _{1}^{\alpha },\kappa _{2}^{-1}\kappa _{2}^{\alpha
},\kappa _{2}^{-1}\kappa _{2}^{\alpha }\right) , \\
\kappa ^{-1}\kappa ^{\beta } &=&\left( \kappa _{1}^{-1}\kappa _{2}^{\beta
},\kappa _{1}^{-1}\kappa _{2}^{\beta },\kappa _{2}^{-1}\kappa _{1},\kappa
_{2}^{-1}\kappa _{1}\right) , \\
\alpha ^{2} &=&\left( \alpha ^{2}z,\alpha ^{2}z,\alpha ^{2},\alpha
^{2}\right) \text{.}
\end{eqnarray*}%
Thus, 
\begin{eqnarray*}
\kappa _{1}^{-1}\kappa _{1}^{\alpha } &=&\kappa _{2}^{-1}\kappa _{2}^{\alpha
}=e, \\
\kappa _{1}^{-1}\kappa _{2}^{\beta } &=&\alpha ^{2}z,\kappa _{2}^{-1}\kappa
_{1}=\alpha ^{2},
\end{eqnarray*}%
\begin{equation*}
\kappa _{1}=\kappa _{2}\alpha ^{2},\kappa _{2}^{-1}\kappa _{2}^{\beta
}=\alpha ^{4}z
\end{equation*}%
and we find that $\kappa _{2}=\alpha \kappa ^{2}$ is a solution. Thus 
\begin{equation*}
\kappa =\left( \alpha ^{3}\kappa ^{2},\alpha ^{3}\kappa ^{2},\alpha \kappa
^{2},\alpha \kappa ^{2}\right)
\end{equation*}%
satisfies all our conditions.

The group $S=\left\langle \alpha ,\beta ,\kappa \right\rangle $ is a $%
\mathfrak{T}_{3}$ group with the presentation%
\begin{equation*}
\left\{ \alpha ,\beta ,\kappa |\left[ \alpha ,\beta ,\alpha \right] =\left[
\alpha ,\beta ,\beta \right] =\left[ \alpha ,\beta ,\kappa \right] =\left[
\alpha ,\kappa \right] =e,\left[ \beta ,\kappa \right] =\alpha ^{-2}\right\} 
\text{.}
\end{equation*}%
The groups $R,S$ are the first two terms of an infinite sequence of
nilpotent subgroups of $C_{\mathcal{A}}\left( z\right) $, which we will
construct in Section 5.4.

\subsection{State-closed abelian groups}

Given $\alpha \in \mathcal{A}$ we indicate the diagonal automorphism $\left(
\alpha ,\alpha ,...,\alpha \right) $ by $\alpha ^{\left( 1\right) }$ and
inductively, $\left( \alpha ^{\left( i\right) },\alpha ^{\left( i\right)
},...,\alpha ^{\left( i\right) }\right) $ by $\alpha ^{\left( i+1\right) }$.

The following theorem shows that recurrent abelian groups (no conditions on
type) are in a sense small.

\begin{theorem}
Let $Y=\left\{ 1,2,...,m\right\} $, $\mathcal{A}=Aut(\mathcal{T(}Y))$.%
\newline
(i) Let $G$ be an abelian recurrent subgroup of $\mathcal{A}$ and $C_{%
\mathcal{A}}\left( G\right) $ be the centralizer of $G$ in $\mathcal{A}$.
Let $\widehat{G}$ the topological closure of $G$ in $\mathcal{A}$. Then, $C_{%
\mathcal{A}}\left( G\right) =\widehat{G}$ .\newline
(ii) Let $m$ be a prime number and $G$ an infinite abelian state-closed
subgroup of $\mathcal{A}$, which acts transitively on the $1$st level of the
tree. Then, $C_{\mathcal{A}}\left( G\right) =\widehat{G}$.
\end{theorem}

\begin{proof}
(i) Let the vertices of $\mathcal{T(}Y)$ be indexed by sequences from $%
Y=\left\{ 1,2,...,m\right\} $. Let $G$ induce the permutation group $P$ on
the set $Y$. Then, $P$ is an abelian transitive permutation group of degree $%
m$ and is therefore regular; it follows that the stabilizer in $G$ of any $%
y\in Y$ is the same as the stabilizer of the $1$st level of the tree $%
H=stab_{G}(1)$. Since the representation of $G$ is recurrent, the projection 
$\pi _{v}$ of $stab_{G}(k)$ on any of its coordinates $v$ produces the group 
$G$.

For every $\sigma \in P$, choose $a_{\mathbf{0}}\left( \sigma \right)
=\left( a_{\mathbf{0}}\left( \sigma \right) _{1},...,a_{\mathbf{0}}\left(
\sigma \right) _{m}\right) \sigma \in G$ which induces $\sigma $ on $Y$. Let 
$h=$ $\left( h_{1},h_{2},...,h_{m}\right) \in $ $H$. Then 
\begin{equation*}
h^{a_{\mathbf{0}}\left( \sigma \right) }=\left( \left( h_{1}\right) ^{a_{%
\mathbf{0}}\left( \sigma \right) _{1}},\left( h_{2}\right) ^{a_{\mathbf{0}%
}\left( \sigma \right) _{2}},...,\left( h_{m}\right) ^{a_{\mathbf{0}}\left(
\sigma \right) _{m}}\right) ^{\sigma }
\end{equation*}%
\begin{equation*}
=\left( h_{1},h_{2},...,h_{m}\right) ^{\sigma }
\end{equation*}%
since $h_{i},a_{\mathbf{0}}\left( \sigma \right) _{i}\in G$ which is
abelian. On varying $\sigma \in P$ we find that $h=$ $\left(
h_{1},h_{1},...,h_{1}\right) $.

Now, for every $\sigma \in P$, there exists $a_{\mathbf{1}}\left( \sigma
\right) =\left( a_{\mathbf{0}}\left( \sigma \right) ,...,a_{\mathbf{0}%
}\left( \sigma \right) \right) \in H$, which induces $\sigma ^{\left(
1\right) }$ modulo $stab_{\mathcal{A}}(2)$. Thus, we produce a sequence $a_{%
\mathbf{i}}\left( \sigma \right) \in stab_{G}(i)$ of elements in $G$ such
that $a_{\mathbf{i}}\left( \sigma \right) =\sigma ^{\left( i\right) }$
modulo $stab_{\mathcal{A}}(i+1)$.

Let $\gamma \in C=C_{\mathcal{A}}\left( G\right) $. Then, 
\begin{eqnarray*}
\gamma &=&\left( \gamma _{1},...,\gamma _{m}\right) \sigma , \\
\gamma ^{\prime } &=&\gamma .a_{\mathbf{0}}\left( \sigma \right)
^{-1}=\left( \gamma _{1}^{\prime },...,\gamma _{m}^{\prime }\right) \in
stab_{C}(1)
\end{eqnarray*}%
and $\gamma _{1}^{\prime }=...=\gamma _{m}^{\prime }$; say $\gamma
_{1}^{\prime }$ induces a permutation $\sigma ^{\prime }$ on $Y.$ Thus, 
\begin{equation*}
\gamma .a_{\mathbf{0}}\left( \sigma \right) ^{-1}.a_{\mathbf{1}}\left(
\sigma ^{\prime }\right) ^{-1}\in stab_{C}(2)\text{.}
\end{equation*}%
We produce in this manner a sequence 
\begin{equation*}
a_{\mathbf{0}}\left( \sigma \right) ,a_{\mathbf{1}}\left( \sigma ^{\prime
}\right) ,a_{\mathbf{2}}\left( \sigma ^{\prime \prime }\right) ,...
\end{equation*}%
of elements of $G$ such that $\gamma $ is equal to the infinite product 
\begin{equation*}
a_{\mathbf{0}}\left( \sigma \right) a_{\mathbf{1}}\left( \sigma ^{\prime
}\right) a_{\mathbf{2}}\left( \sigma ^{\prime \prime }\right) ...\text{.}
\end{equation*}%
Hence, $C_{\mathcal{A}}\left( G\right) =\widehat{G}$.

(ii) Let $m=p$, a prime number. The permutation group $P$ induced on $%
Y=\left\{ 1,...,p\right\} $ is cyclic, say generated by $\sigma $. Since $G$
is infinite, there exists an $h=$ $\left( h_{1},h_{1},...,h_{1}\right) \in $ 
$H$ such that $h_{1}\not\in H$ and therefore we may assume $h_{1}$ induces $%
\sigma $ on $Y$. We produce elements $a_{\mathbf{i}}\in $ $G$ such that $a_{%
\mathbf{i}}=\sigma ^{\left( i\right) }$ modulo $stab_{\mathcal{A}}(i+1)$ and
the proof continues as previously.
\end{proof}

The group $G$ in part (ii) need not be recurrent. For example, let $%
Y=\left\{ 1,2\right\} $, $G$ be the cyclic subgroup of $Aut\left( \mathcal{T(%
}Y)\right) $ generated by $\alpha =\left( e,\alpha ^{3}\right) \sigma \ $%
where $\sigma $ is the transposition $\left( 1,2\right) $.

\begin{proposition}
Let $G$ be a finitely generated abelian group, $H$ a normal subgroup of $G$
such that $\frac{G}{H}$ is cyclic of prime power order. Suppose $(G,H,f)$ is
a simple triple. Then, either $G$ is finite or free abelian.
\end{proposition}

\begin{proof}
Let $Tor\left( G\right) $ be the torsion subgroup of $G$. Then, we have the
decompositions%
\begin{equation*}
G=Tor\left( G\right) \oplus K,H=Tor\left( H\right) \oplus M
\end{equation*}%
where $Tor\left( H\right) \leq Tor\left( G\right) $ and we may assume $M\leq
K$. Suppose $G$ is a mixed group. Then the first possibility is $K=M$ and $%
\frac{Tor\left( G\right) }{Tor\left( H\right) }$ cyclic. Let $n$ be the
exponent of $Tor\left( G\right) $. Then, we have $\left( M^{n}\right)
^{f}=\left( M^{f}\right) ^{n}\leq G^{n}=M^{n}$; a contradiction. The other
possibility is $Tor\left( G\right) =Tor\left( H\right) $ and $\frac{K}{M}$
finite cyclic ; but as $Tor\left( H\right) $ is $f$-invariant, $Tor\left(
G\right) =Tor\left( H\right) =\left\{ e\right\} $ and again we have a
contradiction.
\end{proof}

\section{Sub-triples and Quotient triples}

A subgroup $K$ of $G$ is \textit{semi-invariant} under the action of $f$
provided $\left( K\cap H\right) ^{f}\leq K$. If $K\leq H$ and $K^{f}\leq K$
then $K$ is $f$-\textit{invariant}. Given a triple $\left( G,H,f\right) $
and $G_{1}\leq G,H_{1}\leq H\cap G_{1}$ such that $\left( H_{1}\right)
^{f}\leq G_{1}$, we call $\left( G_{1},H_{1},f|_{H_{1}}\right) $ a \textit{%
sub-triple. }If $N$ is a normal semi-invariant subgroup of $G$ then $%
\overline{f}:\frac{HN}{N}\rightarrow \frac{G}{N}$ given by $\overline{f}%
:Nh\rightarrow Nh^{f}$ is well-defined and $\left( \frac{G}{N},\frac{HN}{N},%
\overline{f}\right) $ is a \textit{quotient triple}.

Given a triple $(G,H,f)$, we produce a sequence of subtriples $%
(G(i),H(i),f_{i})$ defined as follows: 
\begin{equation*}
G(0)=G\text{, }H(0)=H\text{, }f_{0}=f\text{,}
\end{equation*}%
and for $i\geq 1$%
\begin{equation*}
G(i)=H(i-1)^{f}\text{, }H(i)=H(i-1)\cap G(i)\text{, }f_{i}=f_{i-1}|_{H(i)}%
\text{.}
\end{equation*}%
Clearly, if $f$ is an epimorphism, then the sequence stops at $i=0$.

\textbf{Example 1. }The following group $G$ of automorphisms of the binary
tree provides an example for which the sequence $G(i)$\ is infinite,%
\begin{equation*}
G=<\alpha =\left( 1,\alpha \beta ^{2}\right) \sigma \text{, }\beta =\left(
\alpha ,\alpha \right) >\text{.}
\end{equation*}%
It is straightforward to check that $G$ is free abelian of rank $2$, the
subgroup $H=<$ $\alpha ^{2},\beta >$ is of index $2$ and the projection of $%
H $ on the second coordinate is an extension of 
\begin{equation*}
f:\alpha ^{2}\rightarrow \alpha \beta ^{2},\beta \rightarrow \alpha \text{.}
\end{equation*}%
We claim that for $i\geq 1$%
\begin{equation*}
G(i)=<\alpha ^{2^{i-1}},\alpha ^{r_{i}}\beta ^{2}>
\end{equation*}%
where 
\begin{eqnarray*}
r_{1} &=&r_{2}=1, \\
r_{i} &=&1+4t_{i}\text{ such that }t_{i}r_{i-1}\equiv 1\text{ mod }\left(
2^{i-2}\right) \text{ for }i\geq 3\text{.}
\end{eqnarray*}%
The assertion is true for $i=1,2$:

\begin{eqnarray*}
G\left( 0\right) &=&G=<\alpha ,\beta >,H\left( 0\right) =H=<\alpha
^{2},\beta >; \\
G(1) &=&H(0)^{f}=<\alpha \beta ^{2},\alpha >=<\alpha ,\beta ^{2}>, \\
H(1) &=&<\alpha ^{2},\beta ^{2}>; \\
G(2) &=&H(1)^{f}=<\alpha \beta ^{2},\alpha ^{2},>=<\alpha ^{2},\alpha \beta
^{2}>\text{.}
\end{eqnarray*}%
Now, suppose%
\begin{equation*}
G(i)=<\alpha ^{2^{i-1}},\alpha ^{r_{i}}\beta ^{2}>\text{.}
\end{equation*}%
Then 
\begin{eqnarray*}
H(i) &=&<\left( \alpha ^{r_{i}}\beta ^{2}\right) ^{2},\alpha ^{2^{i-1}}>, \\
G(i+1) &=&<\left( \alpha \beta ^{2}\right) ^{r_{i}}\alpha ^{4},\left( \alpha
\beta ^{2}\right) ^{2^{i-2}}>\text{.}
\end{eqnarray*}%
Viewing $G$ as an additive group with basis $\alpha ,\beta $, the generators
of $G(i+1)$ are the rows of the matrix $M=\left( 
\begin{array}{cc}
4+r_{i} & 2r_{i} \\ 
2^{i-2} & 2^{i-1}%
\end{array}%
\right) $.

Let $m,k$ be integers such that $mr_{i}+k2^{i-2}=1$ and let $S=\left( 
\begin{array}{cc}
2^{i-2} & -r_{i} \\ 
m & k%
\end{array}%
\right) $. Then, $\det \left( S\right) =1$ and 
\begin{equation*}
SM=\left( 
\begin{array}{cc}
2^{i} & 0 \\ 
m\left( 4+r_{i}\right) +k2^{i-2} & 2%
\end{array}%
\right) \text{.}
\end{equation*}%
Since $m\left( 4+r_{i}\right) +k2^{i-2}=4m+mr_{i}+k2^{i-2}=4m+1$, we reach 
\begin{equation*}
G(i+1)=<\alpha ^{2^{i}},\alpha ^{1+4m}\beta ^{2}>\text{.}
\end{equation*}

Clearly, $\left[ G:G(i)\right] =2^{i}$ and therefore, $G(i)\not=G(j)$ for $%
i<j$.

\begin{proposition}
Let $G$ be group, $H$ a subgroup $G$ and suppose $(G,H,f)$ is a simple
triple such that $G=Z(G)H^{f}H$. Define, 
\begin{equation*}
G(1)=H^{f}\text{, }H(1)=H\cap H^{f}\text{, }f_{1}=f|_{H(1)}\text{.}
\end{equation*}%
Then, $\left( G(1),H(1),f_{1}\right) $ is a simple triple.
\end{proposition}

\begin{proof}
Let $Y$ be a right transversal of $H$ in $G$ such that $Y$ is contained in $%
Z(G)H^{f}$. Let $K$ be a subgroup of $H(1)$, normal in $G(1)$, with $%
K^{f}\leq K$. Then, $K\leq K^{f^{-1}}$ and $K^{f^{-1}}$ is normal in $H$.
Then, 
\begin{equation*}
K=K^{y}\leq \left( K^{f^{-1}}\right) ^{y}=\left( K^{f^{-1}}\right) ^{hy}
\end{equation*}
for all $y\in Y$ and all $h\in H$; that is, $K\leq \left( K^{f^{-1}}\right)
^{g}$ for all $g\in G$ \ Let $M=\cap _{g\in G}\left( K^{f^{-1}}\right) ^{g}$%
. Then, $M\leq K^{f^{-1}}\leq H$, and $M$ is normal in $G$ such that 
\begin{equation*}
K\leq M\leq K^{f^{-1}}\text{, }M^{f}\leq K\leq M\text{;}
\end{equation*}%
therefore, $M=\left\{ e\right\} =K$.
\end{proof}

The above result generalizes Lemma 3.2 in \cite{NekSid}.

\begin{lemma}
Let $(G,H,f)$ be a triple. Suppose $B$ is semi-invariant and $\sqrt[G]{B}$
is a group. Then $\sqrt[G]{B}$ is semi-invariant.
\end{lemma}

\begin{proof}
Let $x\in \sqrt[H]{B}$, then $x\in H$ and $x^{n}\in B\cap H$. for some $n$.
Therefore, 
\begin{equation*}
\left( x^{n}\right) ^{f}=\left( x^{f}\right) ^{n}\in \left( B\cap H\right)
^{f}
\end{equation*}
and so, $x^{f}\in \sqrt[G]{\left( B\cap H\right) ^{f}}\leq \sqrt[G]{B}$.
\end{proof}

\subsection{Facts about nilpotent groups}

We list below some facts about nilpotent groups which are either well-known
(see \cite{Hall},\cite{Segal}) or have direct proofs. Let 
\begin{equation*}
\left\{ Z_{i}\left( G\right) |1\leq i\leq c\right\} \text{, }\left\{ \gamma
_{i}\left( G\right) |1\leq i\leq c\right\}
\end{equation*}
denote the upper and lower central series of $G$, respectively.

I. Let $G$ be a nilpotent group of class $c$.

1. For all $1\leq i,j\leq c$ 
\begin{eqnarray*}
\left[ \gamma _{i}\left( G\right) ,\gamma _{j}\left( G\right) \right] &\leq
&\gamma _{i+j}\left( G\right) . \\
\left[ Z_{i}\left( G\right) ,\gamma _{j}\left( G\right) \right] &\leq
&Z_{i-j}\left( G\right) \text{.}
\end{eqnarray*}

2. If $G=Z_{i}\left( G\right) H$ for some $i$, then for all $1\leq j\leq c$, 
\begin{eqnarray*}
\gamma _{j}\left( G\right) &\leq &Z_{i-j+1}\left( G\right) \gamma _{j}\left(
H\right) , \\
\gamma _{i+1}\left( H\right) &=&\gamma _{i+1}\left( G\right) \text{.}
\end{eqnarray*}

3. The subset $Tor(G)$ of $G$ of elements of finite order is a subgroup of $%
G $.

4. If $Z(G)\leq Tor(G)$ then $G=Tor(G)$.

5. Suppose $N$ is a normal torsion-free subgroup of $G$. Let $x\in G$, $y\in
N$ and $n$ a positive integer. Then 
\begin{equation*}
\lbrack x^{n},y]=e\Rightarrow \lbrack x,y]=e\text{.}
\end{equation*}

6. Suppose $Tor(G)$ has finite exponent $s$. If $G=Tor(G)K$ for some $K\leq
G $. Then, $G^{s}=K^{s}$.

II. Let $G$ be torsion-free nilpotent.

1. Let $K$ be a subgroup of $G$. Then, the isolator of $K$ in $G$ 
\begin{equation*}
\sqrt[G]{K}=\left\{ x\in G|x^{n}\in K\text{ for some positive integer }%
n\right\}
\end{equation*}%
is a subgroup of $G$. If furthermore $G$ is finitely generated then $\left[ 
\sqrt[G]{K}:K\right] $ is finite.

2. Let $H$ be a subgroup of finite index $m$ in $G$. Then, $H\cap
Z_{i}(G)=Z_{i}(H)$ for all $i$. Also, $\left[ Z_{i}(G):Z_{i}(H)\right]
=q_{i} $ is finite for all $i$ and $q_{i}$ divides $q_{j}$ for $i\leq j$.

III. Let $G$ be finitely generated nilpotent group of class $c$.

1. Then $G$ is Hopfian and has a finite Hirsch length denoted by $h(G)$.
Also, a subgroup $H$ has finite index in $G$ if and only if $h(H)=h(G)$.

2. Let $H$ be subgroup of finite index in $G$. Then $Z(G)$ and $Z(H)$ have
the same Hirsch length. Also, $\left[ \gamma _{i}\left( G\right) :\gamma
_{i}\left( H\right) \right] $ is finite.

\subsection{Triples for nilpotent groups}

\begin{lemma}
Let $G$ be a nilpotent group, $(G,H,f)$ a triple, $\overline{G}=\frac{G}{%
Tor\left( G\right) }$ and $\overline{H}=\frac{HTor\left( G\right) }{%
Tor\left( G\right) }$. Then, $(Tor\left( G\right) ,Tor\left( H\right) ,f)$, $%
(\overline{G},\overline{H},\overline{f})$ are triples. Furthermore, if $G$
is finitely generated and $(G,H,f)$ is simple then $(\overline{G},\overline{H%
},\overline{f})$ is simple.
\end{lemma}

\begin{proof}
The first assertion follows from $Tor\left( G\right) \cap H=Tor\left(
H\right) $. Let $\frac{L}{Tor\left( G\right) }\leq \overline{H}$ be the $%
\overline{f}$-core$\left( \overline{H}\right) $. Then, 
\begin{equation*}
L=Tor\left( G\right) \left( L\cap H\right) =Tor\left( L\right) \left( L\cap
H\right)
\end{equation*}%
and $\left( L\cap H\right) ^{f}\leq L$. By Subsection 3.1, item I.6, there
exists $s\geq 1$ such that $L^{s}=\left( L\cap H\right) ^{s}$. Therefore, $%
L^{s}$ is $f$-invariant. As $f$ is simple, we have $L^{s}=\left\{ e\right\} $%
, $L=Tor\left( G\right) $.
\end{proof}

\begin{lemma}
Let$\ G$ be a $\mathfrak{T}$ -group, $H$ a subgroup of $G$ of finite index
in $G$ and $f:H\rightarrow G\ $a monomorphism. Then $H^{f}$ has finite index
in $G$. If $U$ is an $f$-invariant normal subgroup of $H$ then $U\cap Z(H)$
is an $f$-invariant normal subgroup of $G$. If $\left( G,H,f\right) $ is
simple then $\left( \left\langle H,H^{f}\right\rangle ,H,f\right) $ is
simple. .
\end{lemma}

\begin{proof}
The groups $H$ and $G$ have equal Hirsch lengths, by Subsection 3.1, item
I.4. Since $H^{f}\cong H$, it has the same Hirsch length as $G$ and
therefore $H^{f}$ has finite index in $G$. Therefore, $Z(H),Z(H)^{f}\leq
Z(G) $. Let $W=U\cap Z(H)$. Then, 
\begin{equation*}
W^{f}\leq U^{f}\cap Z(H)^{f}\leq U\cap Z(G)=W\text{.}
\end{equation*}%
The last assertion follows directly.
\end{proof}

\begin{proposition}
Let $G$ be a $\mathfrak{T}_{c}$ -group, $H$ a subgroup of finite index in $G$
, $f\in Hom\left( H,G\right) $ be a simple monomorphism and $L\leq G$ be
defined by $\frac{L}{Z(G)}=\overline{f}$-core$(\frac{HZ(G)}{Z(G)})$. Then 
\newline
(i) $\overline{f}:\frac{HZ(G)}{Z(G)}\rightarrow \frac{G}{Z(G)}$ induced by $%
f $ is a monomorphism;\newline
(ii) $L,\sqrt[G]{L}$ are abelian, semi-invariant and the corresponding
quotient triples\newline
$\left( \frac{G}{L},\frac{HL}{L},\overline{f}\right) $,$\left( \frac{G}{\sqrt%
[G]{L}},\frac{H\sqrt[G]{L}}{\sqrt[G]{L}},\overline{f}\right) $ are simple;%
\newline
(iii) if $HL=H\sqrt[G]{L}$ then $L=\sqrt[G]{L}$;\newline
(iv) if $G=H\sqrt[G]{L}$ then $G$ is abelian.
\end{proposition}

\begin{proof}
Let $h\in H$ such that $Z(G)h\in \ker \left( \overline{f}\right) $. Then, $%
h^{f}\in Z(G)$. As $Z(H)^{f}=Z\left( H^{f}\right) \leq Z(G)$ and $\ker
\left( f\right) =\left\{ e\right\} $, it follows that $h\in Z(H)\leq Z(G)$.

We have $Z(G)\leq L\leq Z(G)H$. Let $M=L\cap H$. Then, $M$ is a normal
subgroup of $H$ and $L=Z(G)+M,L^{\prime }=M^{\prime }$. Also, $\left(
Z(G)x\right) ^{\overline{f}}=Z(G)x^{f}\in Z(G)M$ for all $x\in M$; that is, $%
M^{f}\leq Z(G)M$. Therefore $\left( M^{\prime }\right) ^{f}\leq M^{\prime }$%
. Since $f$ is simple, $M$ and $L$ are abelian and therefore $\sqrt[G]{L}$
is abelian.

Write $\overline{G}=\frac{G}{Z(G)},\overline{L}=\frac{L}{Z(G)}$. As $\sqrt[G]%
{L}$ is abelian, easily, $\sqrt[\overline{G}]{\overline{L}}=\frac{\sqrt[G]{L}%
}{Z(G)}$. If $x\in \sqrt[H]{L}$ then there exists $n$ such that $x^{n}\in
L\cap H=M$; therefore $\left( x^{n}\right) ^{f}=\left( x^{f}\right) ^{n}\in
L $ and $x^{f}\in \sqrt[G]{L}$. The assertion that $\left( \frac{G}{\sqrt[G]{%
L}},\frac{H\sqrt[G]{L}}{\sqrt[G]{L}},\overline{f}\right) $ is simple is now
clear.

Let $\overline{G}=\frac{G}{L}$and $\overline{H}=\frac{HL}{L}$. Then $\frac{%
\sqrt[G]{L}}{L}=Tor\left( \overline{G}\right) $ and it follows that from
Lemma 2 that $\left( \frac{G}{\sqrt[G]{L}},\frac{H\sqrt[G]{L}}{\sqrt[G]{L}},%
\overline{f}\right) $ is simple. Suppose $HL=H\sqrt[G]{L}$. Then, $\overline{%
H}=\overline{H}Tor\left( \overline{G}\right) $ and the equalities $Tor\left( 
\overline{H}\right) =Tor\left( \overline{G}\right) =\left\{ L\right\} $
follow; that is, $L=\sqrt[G]{L}$.

Suppose $G=\sqrt[G]{L}H$. Then working modulo $L$, we have $\sqrt[G]{L}%
=T\left( G\right) $ and therefore there exists$\ s\geq 1$ such that $%
G^{s}=H^{s}$. In other words, going back to $G$, we have $G^{s}L=H^{s}L$. It
follows from $L=Z(G)M$ that 
\begin{equation*}
\left[ L,G^{s}\right] =\left[ L,H^{s}\right] =\left[ M,G^{s}\right] =\left[
M,H^{s}\right] \text{;}
\end{equation*}%
hence, $\left[ M,H^{s}\right] $ is an $f$-invariant subgroup. Therefore,%
\begin{equation*}
\left[ M,H^{s}\right] =\left[ M,G^{s}\right] =\left\{ e\right\} ,M\leq Z(G)
\end{equation*}%
and 
\begin{equation*}
L=Z(G)=\sqrt[G]{L},G=Z(G)H,G^{\prime }=H^{\prime }=\left\{ e\right\} \text{.}
\end{equation*}
\end{proof}

We show in the next example that$\sqrt[G]{L}$ can be different from $L$.

\textbf{Example 2.} Let $G=F\left( 2,2\right) $ freely generated by $%
x_{1},x_{2}$. Let $n>1$ and let $H=<x_{1}^{n},x_{2}^{n}>$. Then $Z\left(
H\right) =\left\langle \left[ x_{2},x_{1}\right] ^{n^{2}}\right\rangle $.
The map $f:x_{1}^{n}\rightarrow x_{1}^{n},x_{2}^{n}\rightarrow x_{2}$
extends to a monomorphism from $H$ into $G$ where $f:\left[ x_{2},x_{1}%
\right] ^{n^{2}}\rightarrow \left[ x_{2},x_{1}\right] ^{n}$. It is clear
then that $f$ is simple, $L=Z(G)<x_{1}^{n}>$ and $\sqrt[G]{L}=Z(G)<x_{1}>$.

\section{Kernel versus Core}

\begin{proposition}
Let $K,P$ be groups, $P$ a transitive permutation group on the set $%
Y=\left\{ 1,...,m\right\} $ and $P_{1}$ be the stabilizer of $1$ in $P$.
Furthermore, let $W$ be the wreath product $Kwr_{Y}P$ and let $W$ act on $Y$
as $P$. Let $B=K^{Y}$ and $W_{1}=$ $BP_{1}$. Consider a nilpotent subgroup $%
G $ of $W$ which induces a transitive group on $Y$. Let $x=\left(
x_{1},x_{2},...,x_{m}\right) \sigma \in G_{1}=G\cap W_{1}$. If $x_{1}$ has
finite order then $x$ also has finite order.
\end{proposition}

\begin{proof}
Suppose that there exists an $x=\left( x_{1},x_{2},...,x_{m}\right) \sigma $
where $\sigma \in P_{1}$ such that $x_{1}$ has finite order, yet $o\left(
x\right) $ is infinite. Then $x^{\prime }=\left( x^{o\left( x_{1}\right)
}\right) ^{o\left( \sigma \right) }\in B_{G}=G\cap B$ and $x^{\prime
}=\left( e,x_{2}^{\prime },...,x_{m}^{\prime }\right) $ has infinite order.
Let $X$ be the set of non-trivial $x\in $ $B_{G}$ such that each $x_{i}=e$
or $o\left( x_{i}\right) $ is infinite. Choose $x\in X$ such that first, $x$
has a maximum number of trivial entries and second, $x\in \gamma _{j}\left(
G\right) $ for a maximum $j$. Let $g=\left( g_{1},g_{2},...,g_{m}\right)
\rho \in G$ where $\rho \in P$ and $o\left( \rho \right) =r$. Then, $\left[
x,g^{r}\right] \in B_{G}$ has at least the same number of trivial entries as 
$x$ and $\left[ x,g^{r}\right] \in \gamma _{j+1}\left( G\right) $.
Therefore, $\left[ x,g^{r}\right] $ has finite order. Hence, in the
torsion-free nilpotent group $\overline{G}=\frac{G}{Tor\left( G\right) }$,
we have $\left[ \overline{x},\overline{g^{r}}\right] =\overline{e}$ and so, $%
\left[ \overline{x},\overline{g}\right] =\overline{e}$; that is, $\left[ x,g%
\right] \ $has finite order. We may assume $x_{1}=e$ and let $i$ be such $%
x_{i}$ has infinite order. Now let $g=\left( g_{1},g_{2},...,g_{m}\right)
\rho \in G$ be such that $\left( i\right) \rho =1$. Then,%
\begin{eqnarray*}
\left[ x,g\right] &=&x^{-1}x^{g}=\left( e,x_{2}^{-1},...,x_{m}^{-1}\right)
\left( e,x_{2}^{g_{2}},...,x_{m}^{g_{m}}\right) ^{\rho } \\
&=&\left( x_{i}^{g_{i}},\ast ,...,\ast \right)
\end{eqnarray*}%
which has infinite order; a contradiction is reached.
\end{proof}

\begin{theorem}
Let $G$ be a nilpotent group, $H$ a subgroup of finite index $m$ in $G$, $%
f\in Hom\left( H,G\right) $ and $L=$ $f$-$core(H)$. Then,\newline
(i) $\ker (f)\leq \sqrt[H]{L}$;\newline
(ii) $Tor\left( H\right) ^{f}\leq Tor\left( H^{f}\right) \leq \left( \sqrt[H]%
{L}\right) ^{f}$; \newline
(iii) if $L=\left\{ e\right\} $ then $Tor\left( H\right) ^{f}=Tor\left(
H^{f}\right) $; \newline
(iv) if $G$ is finitely generated and $f$ an epimorphism, then $\sqrt[G]{L}%
=L $.
\end{theorem}

\begin{proof}
The triple $(G,H,f)$ provides us with a state-closed representation $\varphi
:\frac{G}{L}\rightarrow Aut\left( \mathcal{T}\left( Y\right) \right) $, for $%
Y=\left\{ 1,...,m\right\} $ and where for $h\in H$, we have 
\begin{equation*}
h^{\varphi }=\left( h^{f\varphi },\ast ,...,\ast \right) \sigma
\end{equation*}%
and $\left( 1\right) ^{\sigma }=1$.

(i) If $h\in \ker \left( f\right) $ then $h^{\varphi }=\left( e,\ast
,...,\ast \right) \sigma $ and by the previous proposition, $h^{\varphi }$
has finite order. As $L=\ker \left( \varphi \right) $, we have $h\in \sqrt[H]%
{L}$ and we are done.

(ii) The first inclusion is clear. Now suppose $x\in Tor\left( H^{f}\right) $%
; that is, $x=h^{f}$ and $x^{n}=e$ for some $n$. Then, 
\begin{eqnarray*}
e &=&(h^{f})^{n}=(h^{n})^{f},h^{n}\in \ker \left( f\right) , \\
h &\in &\sqrt[H]{\ker \left( f\right) }\leq \sqrt[H]{L}\text{, }x\in (\sqrt[H%
]{L})^{f}\text{.}
\end{eqnarray*}%
If $L=\left\{ e\right\} $ then $\sqrt[H]{L}=Tor\left( H\right) $ and the
result follows from $Tor\left( H\right) ^{f}\leq Tor\left( H^{f}\right) \leq
Tor\left( H\right) ^{f}$.

(iii) follows immediately from (ii).

(iv) Since $G$ is finitely generated, $Tor(G)$ is a finite group. Suppose
initially that $L$ is trivial. Then, $\sqrt[G]{L}=Tor(G)$ and by item (ii), $%
Tor\left( H\right) ^{f}\leq Tor\left( G\right) \leq Tor\left( H\right) ^{f}$%
; thus,

$Tor\left( H\right) ^{f}=Tor(G)$. As, $Tor\left( H\right) \leq Tor(G)$, we
have $Tor\left( H\right) =Tor(G)$ and as $f$ is simple, we conclude that $%
Tor(G)$ is trivial.

In the general case, we consider the triple $(\frac{G}{L},\frac{H}{L},%
\overline{f})$. Then, $\overline{f}$ is a simple epimorphism and therefore $%
Tor\left( \frac{G}{L}\right) =\left\{ L\right\} $; that is, $\sqrt[G]{L}=L$.
\end{proof}

\begin{corollary}
Let $G$ be a torsion-free nilpotent group, $H$ a subgroup of finite index $m$
and $f$:$H\rightarrow G$ a homomorphism. Then%
\begin{equation*}
\text{ }f\text{ simple}\Rightarrow \ker \left( f\right) =\left\{ e\right\} 
\text{.}
\end{equation*}%
Suppose $\ker \left( f\right) =\left\{ e\right\} $. Then 
\begin{equation*}
f\text{ simple}\Leftrightarrow f:Z(H)\rightarrow Z(G)\text{ simple.}
\end{equation*}
\end{corollary}

\begin{proof}
The first assertion is a direct application of part (i) of the theorem.
Suppose $G$ is finitely generated; then $Z(H)\leq Z(G)$. It follows easily
that $f$ simple implies that $f|_{Z(H)}:Z(H)\rightarrow Z(G)$ simple. On the
other hand, suppose $f:Z(H)\rightarrow Z(G)$ is simple and let $K$ be a
nontrivial subgroup of $H$, normal in $G$ and $f$-invariant. Then, $K\cap
Z(H)=\left\{ e\right\} $ and 
\begin{eqnarray*}
\left( K\cap Z(H)\right) ^{f} &\leq &K\cap Z(H)^{f}=K\cap Z\left(
H^{f}\right) \\
&\leq &K\cap Z\left( G\right) =K\cap H\cap Z\left( G\right) =K\cap Z(H)\text{%
.}
\end{eqnarray*}
\end{proof}

\textbf{Examples 3. }

(1 ) \textit{A simple triple }$\left( G,H,f\right) $\textit{\ where }$G$%
\textit{\ is finite and }$\ker \left( f\right) \not=\left\{ e\right\} $.

Let $p$ be a prime number, $Y=\left\{ 1,2,...,p\right\} $ and $\sigma $ the
permutation $\left( 1,2,...,p\right) $. Let $W_{s}$ be the group of
automorphisms of the $p$-adic tree $\mathcal{T}\left( Y\right) $ generated by

\begin{equation*}
\sigma _{0}=\sigma ,\sigma _{1}=\left( e,...,e,\sigma _{0}\right)
,...,\sigma _{s}=\left( e,...,e,\sigma _{s-1}\right) \text{,}
\end{equation*}%
The $W_{s}$ is the $s$-iterated wreath product $\left(
((C_{p}wr..)wr)C_{p}\right) wrC_{p}(=W_{s-1}wrC_{p})$. Let $%
H=stab_{W_{s}}(1) $ and $\pi _{1}:H\rightarrow W_{s}$. Then, $\left[ W_{s}:H%
\right] =p$, $\pi _{1}\left( H\right) \cong W_{s-1}$ and $\ker \left(
f\right) =\left\{ e\right\} \times W_{s-1}\times ...\times W_{s-1}$.

(2) \textit{A simple triple }$\left( G,H,f\right) $\textit{\ where }$G$%
\textit{\ is of mixed type and }$\ker \left( f\right) \not=\left\{ e\right\} 
$.

Let $G=\left( CwrD\right) \left\langle x\right\rangle $ where $%
C=\left\langle c\right\rangle ,D=\left\langle d\right\rangle $ each of order 
$p$, and $x$ of infinite order inducing conjugation by $d$ on $CwrD$.
Therefore, $G$ is a nilpotent group with $Z\left( G\right) =\left\langle
z,dx^{-1}\right\rangle $ where $z=cc^{d}..c^{d^{p-1}}$. Let $H=\left\langle
C^{D},x^{p}\right\rangle =\left\langle
c,c^{d},..,c^{d^{p-2}},z,x^{p}\right\rangle $ and $M=\left\langle
d,x\right\rangle $ an abelian group of type $\mathbb{Z}_{p}\times \mathbb{Z}$%
. Then, $H$ is abelian of type $\left( \mathbb{Z}_{p}\right) ^{p}\times 
\mathbb{Z}$, $\left[ G:H\right] =p^{2},$ $Z\left( H\right) =\left\langle
z,x^{p}\right\rangle $. The extension of the map%
\begin{equation*}
c\rightarrow 1,c^{d}\rightarrow 1,...,c^{d^{p-2}}\rightarrow 1,z\rightarrow
d,x^{p}\rightarrow x
\end{equation*}
produces an epimorphism $f:H\rightarrow M$. Then, $\ker \left( f\right)
=\left\langle c,c^{d},..,c^{d^{p-2}}\right\rangle $. Note that the only
subgroup of $\ker \left( f\right) $ which is normal in $\left\langle
c,d\right\rangle $ is the trivial subgroup. Let $K$ be an $f$-invariant
subgroup of $H$ normal in $G$. Then, 
\begin{equation*}
K^{f}\leq K\cap M\leq H\cap M=\,\left\langle x^{p}\right\rangle \text{.}
\end{equation*}%
Therefore, $K^{f}=\left\{ e\right\} ,K\leq \ker \left( f\right) $ and so, $%
K=\left\{ e\right\} $.

(3)\textbf{\ }\textit{A triple} $\left( G,H,f\right) $ \textit{where} $G$ 
\textit{is} a $\mathfrak{T}$-\textit{group and} $f$\textit{-core}$\left(
H\right) =Z(H)$.

Let $G=F\left( 2,2\right) $ freely generated by $a,b$. Let $H$ be the
subgroup generated by $a^{3},b^{2}$. Then, $\left[ G:H\right] =36$. Define
the endomorphism $f:H\rightarrow G$ extended from $f:a^{3}\rightarrow
a^{2},b^{2}\rightarrow b^{3}$. Then, $f:\left[ a^{3},b^{2}\right]
\rightarrow \left[ a^{2},b^{3}\right] $. Thus, $f$-core$\left( H\right)
=\left\langle \left[ a,b\right] ^{6}\right\rangle =Z(H)$.

\section{Simple triples for $\mathfrak{T}$-groups}

\subsection{Nilpotent groups: $2$-generated or of class $2$}

\begin{lemma}
Let $G$ be a $\mathfrak{T}_{2}$-group. Then, there exists a subgroup $K$ of $%
G$ such that $G=Z\left( G\right) K,K^{\prime }=Z\left( K\right) $.
\end{lemma}

\begin{proof}
As $Z\left( G\right) $ is isolated in $G$, then $\frac{Z\left( G\right) }{%
G^{\prime }}$ has a complement $\frac{K}{G^{\prime }}$ in $\frac{G}{%
G^{\prime }}$. Therefore,%
\begin{eqnarray*}
G &=&Z\left( G\right) K,G^{\prime }=K^{\prime }\text{,} \\
G^{\prime } &=&Z\left( G\right) \cap K=Z\left( K\right) \text{.}
\end{eqnarray*}
\end{proof}

\begin{theorem}
Let $G$ be a $\mathfrak{T}_{c}$-group with $c\leq 2$ and let $H$ be a
subgroup of finite index in $G$. Then there exists a subgroup $K$ of finite
index in $H$, which admits a strongly simple epimorphism $f:K\rightarrow G$.
\end{theorem}

\begin{proof}
The case $G$ abelian is obvious; so let $G$ have nilpotency class $2$.
Choose $\left\{ Z(G)a_{1},...,Z(G)a_{d}\right\} $, a free generating set of $%
\frac{G}{Z(G)}$ and let $A=\left\langle a_{1},a_{2},...,a_{d}\right\rangle $%
. Then, by the previous lemma, $Z(A)=A^{\prime }$ and $Z(G)=U_{0}\oplus
U_{1} $ where $U_{1}=\sqrt[G]{Z(A)}$. Then, $G=U_{0}\oplus U_{1}A$.

There exists a generating set $\left\{ x_{1},x_{2},...,x_{d}\right\} $ of $A$
such that modulo $Z\left( G\right) $, we have $H=\left\langle
x_{1}^{k_{1}},x_{2}^{k_{2}},...,x_{d}^{k_{d}}\right\rangle $ where $%
k_{i}\geq 1$. Thus, there exist $c_{i}\in Z\left( G\right) $ such that $%
b_{i}=c_{i}x_{i}^{k_{i}}\in H$ $\left( 1\leq i\leq d\right) $. Define $%
B=\left\langle b_{1},b_{2},...,b_{d}\right\rangle $. Then, 
\begin{eqnarray*}
H &=&BZ\left( H\right) , \\
Z\left( B\right) &=&B^{\prime }=\left\langle \left[ x_{i},x_{j}\right]
^{k_{i}k_{j}}|i<j\right\rangle , \\
B^{\prime } &\leq &A^{\prime }=\left\langle \left[ x_{i},x_{j}\right]
|i<j\right\rangle \text{.}
\end{eqnarray*}%
Let $V_{1}=\sqrt[H]{Z(B)}$. Then, $V_{1}\leq U_{1}$ and $\left[ U_{1}:V_{1}%
\right] =r_{1}$ is finite. Now, we prove that we may choose $U_{0}$ such
that $H=V_{0}\oplus V_{1}B$ where $V_{0}\leq U_{0}$. We argue in $\frac{%
Z\left( G\right) }{V_{1}}$. Since $\frac{U_{1}}{V_{1}}=Tor\left( \frac{%
Z\left( G\right) }{V_{1}}\right) $, there exists $W_{0}\leq $ $Z\left(
G\right) $ such that $\frac{Z\left( G\right) }{V_{1}}=\frac{W_{0}\oplus V_{1}%
}{V_{1}}\oplus \frac{U_{1}}{V_{1}}$ and $\frac{Z\left( H\right) }{V_{1}}\leq 
\frac{W_{0}\oplus V_{1}}{V_{1}}$. It follows from $V_{1}\leq Z\left(
H\right) \leq W_{0}\oplus V_{1}$ that $Z\left( H\right) =\left( Z\left(
H\right) \cap W_{0}\right) \oplus V_{1}$. Let $\left[ U_{0}:V_{0}\right]
=r_{0}$.

Now let $r=\func{lcm}\left( r_{0},r_{1}\right) $, $k=\func{lcm}\left\{
k_{i}|1\leq i\leq d\right\} $. Define the subgroups 
\begin{eqnarray*}
B_{0} &=&\left\langle \left( c_{i}x_{i}^{k_{i}}\right) ^{r\frac{k}{k_{i}}%
}|1\leq i\leq d\right\rangle , \\
K &=&Z\left( G\right) ^{r^{2}k^{2}}B_{0}\text{.}
\end{eqnarray*}%
Then $K$ is a subgroup of finite index in $H$ and 
\begin{equation*}
K=U_{0}^{r^{2}k^{2}}\oplus U_{1}^{r^{2}k^{2}}B_{0}\text{.}
\end{equation*}%
Now consider the map%
\begin{equation*}
\gamma :z\rightarrow z^{r^{2}k^{2}},x_{i}\rightarrow \left(
c_{i}x_{i}^{k_{i}}\right) ^{r\frac{k}{k_{i}}}\text{.}
\end{equation*}%
We note that if the map $\gamma $ extends to an endomorphism from $G$ onto $%
K $ then 
\begin{eqnarray*}
\gamma &:&U_{0}\rightarrow U_{0}^{r^{2}k^{2}},U_{1}\rightarrow
U_{1}^{r^{2}k^{2}}, \\
\left[ x_{i},x_{j}\right] &\rightarrow &\left[ x_{i},x_{j}\right]
^{r^{2}k^{2}}\text{ }\left( 1\leq i\leq d\right) \text{.}
\end{eqnarray*}

To prove that $\gamma $ extends to an endomorphism, it is sufficient to
observe that if for some $u_{1}\in U_{1}$ and integer $s$ we have $u_{1}^{s}$
is a word $w\left( \left[ x_{i},x_{j}\right] \right) $ in the commutators $%
\left[ x_{i},x_{j}\right] $ then we have in the extension, $\gamma
:u_{1}^{s}\rightarrow \left( u_{1}^{r^{2}k^{2}}\right) ^{s}=$ $\left(
u_{1}^{s}\right) ^{r^{2}k^{2}}$on the one hand and $\gamma :w\left( \left[
x_{i},x_{j}\right] \right) \rightarrow w\left( \left[ x_{i},x_{j}\right]
^{r^{2}k^{2}}\right) =w\left( \left[ x_{i},x_{j}\right] \right)
^{r^{2}k^{2}} $ on the other and the two images coincide.

\textbf{Example 4. }Given a simple triple $\left( G,H,f\right) $ where $G$
is a $\mathfrak{T}$-group, a question may be posed as to whether assuming $%
f|_{Z(H)}:$ $Z(H)\rightarrow Z(G)$ is an epimorphism implies that $f$ itself
is an epimorphism. The following is a counterexample: let $G=F\left(
2,2\right) $, freely generated by $a,b$. Write $\left[ a,b\right] =z$ and
let $H_{1}=\left\langle a^{4},b^{4},z^{4}\right\rangle $, $%
H_{2}=\left\langle a^{2},b^{2},z\right\rangle $. Then, $f:a^{4}\rightarrow
a^{2},b^{4}\rightarrow b^{2},z^{4}\rightarrow z$ extends to an isomorphism
from $H_{1}$ onto $H_{2}$ and $f$ is simple. Therefore, $Z\left(
H_{1}\right) ^{f}=Z\left( G\right) $, yet $f$ is not an epimorphism.
\end{proof}

\begin{lemma}
Let $G$ is be a $2$-generated $\mathfrak{T}_{c}$-group. Then $G^{\prime
}=Z_{c-1}\left( G\right) $.
\end{lemma}

\begin{proof}
We have $G^{\prime }\leq Z_{c-1}\left( G\right) $ and $\frac{G}{G^{\prime }},%
\frac{G}{Z_{c-1}\left( G\right) }$ are $2$-generated non-cyclic abelian
groups. Since $\frac{G}{Z_{c-1}\left( G\right) }$ is a is torsion-free
quotient of $\frac{G}{G^{\prime }}$, the result follows.
\end{proof}

\begin{theorem}
Let $G$ be $2$-generated $\mathfrak{T}_{c}$-group. Suppose $H$ is a proper
normal subgroup of $G$ of finite index $m$, which is isomorphic to $G$.
Then, $G$ is abelian.
\end{theorem}

\begin{proof}
Let $G$ be generated by $a_{1},a_{2}$ and suppose $c\geq 2$.

First, we will argue the case $c=2$ . Then $G$ is isomorphic to $F\left(
2,2\right) $ and $G^{\prime }=Z(G)=<[a_{1},a_{2}]>$. We may choose the
generators $a_{1},a_{2}$ such that $HZ(G)=<a_{1}^{m_{1}},a_{2}^{m_{2}}>Z(G)$%
. Thus, 
\begin{equation*}
H^{\prime }=Z(H)=<[a_{1},a_{2}]^{m_{1}m_{2}}>\text{.}
\end{equation*}%
We have from $\left[ HZ(G),a_{i}\right] $ $\left( i=1,2\right) $ 
\begin{eqnarray*}
\lbrack a_{1}^{m_{1}},a_{2}] &=&[a_{1},a_{2}]^{m_{1}}\text{, }%
[a_{1},a_{2}^{m_{2}}]=[a_{1},a_{2}]^{m_{2}}\text{,} \\
\lbrack a_{1},a_{2}]^{m_{1}},[a_{1},a_{2}]^{m_{2}} &\in &H\cap Z(G)=Z(H)%
\text{.}
\end{eqnarray*}%
Therefore, $m_{1}m_{2}=\pm 1$. Hence, $HG^{\prime }=G=H$.

Let $c\geq 2$ and $f:H\rightarrow G$ be an epimorphism. Then $f$ induces an
epimorphism $\overline{f}:\frac{HZ_{c-2}(G)}{Z_{c-2}(G)}\rightarrow \frac{G}{%
Z_{c-2}(G)}$, as $H\cap Z_{c-2}(G)=Z_{c-2}(H)$. The class $2$ case leads to $%
G=HZ_{c-2}(G)$. Therefore%
\begin{equation*}
\frac{G}{Z_{c-2}(H)}=\frac{H}{Z_{c-2}(H)}\oplus \frac{Z_{c-2}(G)}{Z_{c-2}(H)}%
\text{.}
\end{equation*}%
Moreover, since $\frac{G}{Z_{c-2}(H)}$ and $\frac{H}{Z_{c-2}(H)}$ are $2$%
-generated, we reach $Z_{c-2}(G)=Z_{c-2}(H)$ and thus, $G=H$; again, a
contradiction.
\end{proof}

\subsection{Derived length and nilpotency class}

\begin{theorem}
Let $G$ be a $\mathfrak{T}$-group, $H$ a subgroup of finite index $m>1$ and $%
f:H\rightarrow G$ simple. Then $s\left( G\right) \leq l\left( m\right) $.
\end{theorem}

\begin{proof}
We may suppose $G$ non-abelian. By Corollary 1, $f$ is a monomorphism. As $%
H^{f}$ is a subgroup of finite index in $G$, we have $Z(H),Z(H)^{f}\leq Z(G)$%
. Clearly, $Z(G)$ is not contained in $H$; for otherwise, $Z(G)=Z(H)$ and $f$%
-invariant.

Consider the triple $\left( \frac{G}{Z(G)},\frac{HZ(G)}{Z(G)},\overline{f}%
\right) $. Then the index $\left[ \frac{G}{Z(G)}:\frac{HZ(G)}{Z(G)}\right] $
is a proper divisor of $m$. Define $L\leq G$ by $\frac{L}{Z(G)}=$ $\overline{%
f}$-core($\frac{HZ(G)}{Z(G)}$). By Proposition 3, both $L,\sqrt[G]{L}$ are
abelian and the triples $\left( \frac{G}{L},\frac{HL}{L},\overline{f}\right)
,\left( \frac{G}{\sqrt[G]{L}},\frac{H\sqrt[G]{L}}{\sqrt[G]{L}},\overline{f}%
\right) $ are simple.

Now, we consider the chain of subgroups 
\begin{equation*}
H\leq HL\leq H\sqrt[G]{L}\leq G\text{. }
\end{equation*}

Since $HZ\left( G\right) \leq $ $HL$ and $Z\left( G\right) \not\leq H$, we
have $H$ is a proper subgroup of $HL$. By Proposition 3, if $HL=H\sqrt[G]{L}$
then $L=\sqrt[G]{L}$ and since $G$ is non-abelian, $H\sqrt[G]{L}\not=G$.

We apply induction on $l(m)$. If $L=\sqrt[G]{L}$ then $\frac{G}{L}$ is
torsion-free, $\left[ \frac{G}{L}:\frac{HL}{L}\right] =m^{\prime }$ and $%
l\left( m^{\prime }\right) <l\left( m\right) $; therefore $s\left( \frac{G}{L%
}\right) \leq $ $l\left( m^{\prime }\right) $ and since $L$ is abelian, $%
s\left( G\right) \leq l\left( m^{\prime }\right) +1\leq l\left( m\right) $.
If $L\not=\sqrt[G]{L}$ then $HL\not=H\sqrt[G]{L}$ and therefore $\left[ 
\frac{G}{\sqrt[G]{L}}:\frac{H\sqrt[G]{L}}{\sqrt[G]{L}}\right] =m^{\prime
\prime }$, $l\left( m^{\prime \prime }\right) <l\left( m\right) $ and the
argument proceeds as in the previous case.
\end{proof}

\begin{corollary}
Let $G$ be a finitely generated nilpotent group. Suppose $\left[ G:H\right]
=p$ a prime number and $\left( G,H,f\right) $ a simple triple. Then, $G$ is
a finite $p$-group (no restriction on nilpotency class or derived length) or
is free abelian.
\end{corollary}

\begin{proof}
Proceed by induction on the order of $Tor(G)$. If $G$ is torsion-free then
by the previous theorem, $G$ is abelian.

If $Tor(G)\neq \left\{ e\right\} $ then $Tor(G)$ is not contained in $H$.
Since $H$ is a maximal subgroup of $G$, by Proposition 2, $\left(
H^{f},H\cap H^{f},f\right) $ is a simple triple of degree $p$. Since $%
\left\vert Tor(H^{f})\right\vert <\left\vert Tor(G)\right\vert $ we conclude
that $H^{f}$ is finite or torsion-free. In the first case, it follows that $%
G $ is finite. In the second case, $\ker (f)=Tor(H)$, normal in $G$ and $f$%
-invariant; therefore, $\ker (f)=\left\{ e\right\} $ and $G$ is
torsion-free, contrary to the assumption.

To justify that neither $s\left( G\right) $ nor $c\left( G\right) $ can be
bounded in case $G$ is finite, we recall $W_{s}$, the iterated wreath
product of cyclic groups of order $p$ in Examples 3 (1).
\end{proof}

\begin{theorem}
Let $G$ be a $\mathfrak{T}$-group, $H$ a subgroup of finite index $m$.
Suppose $f:H\rightarrow G$ is strongly simple. Then $\overline{f}:\frac{HZ(G)%
}{Z\left( G\right) }\rightarrow \frac{G}{Z\left( G\right) }$ is also
strongly simple. Furthermore, $c\left( G\right) \leq l\left( m\right) $.
\end{theorem}

\begin{proof}
Let $G$ be non-ablelian. Suppose there exists a nontrivial subgroup $\frac{L%
}{Z\left( G\right) }$ of $\frac{Z(G)H}{Z\left( G\right) }$ which is $%
\overline{f}$- invariant; then there exists $K\leq H$ such that $L=$ $%
Z\left( G\right) K$; thus $L^{\prime }=K^{\prime }$ $\ $and is $f$%
-invariant. Since $f$ is strongly simple, $L$ is free abelian of finite rank
and we may assume $L=Z(G)\oplus K$, as $Z(G)$ is isolated. The decomposition
of $L$\ provides us with uniquely defined homomorphisms%
\begin{equation*}
\zeta \in Hom\left( K,Z(G)\right) ,\gamma \in Hom\left( K,K\right)
\end{equation*}%
defined by $k^{f}=k^{\zeta }+k^{\gamma }$ for all $k\in K$. We note that $%
\gamma $ is a monomorphism whose characteristic polynomial $\mu $ is
non-constant..

Since $\left[ L:Z(H)\oplus K\right] $ is finite, $f$ extends to an
automorphism $\widehat{f}$ of $\mathbb{Q\otimes }L$, say having
characteristic polynomial $\rho $; clearly, $\mu $ is a factor of $\rho $.
However, since $f$ is simple, $\rho $ does not have monic integral
polynomial factors of positive degree; a contradiction.

Induction on $c\left( G\right) $ leads directly to a proof of the last
assertion.
\end{proof}

The following natural example shows that the above limit is satisfactory.

\textbf{Example 5. }Let $V_{n}$ be the additive free $\mathbb{Z}$-module of
rank $n$ generated by $v_{i}$ $\left( 1\leq i\leq n\right) $ and let $%
x_{n}\in GL\left( V_{n}\right) $ be defined by 
\begin{eqnarray*}
x_{n} &:&v_{i}\rightarrow v_{i}+v_{i+1\text{ }}\text{ }\left( 1\leq i\leq
n-1\right) , \\
v_{n} &\rightarrow &v_{n}\text{.}
\end{eqnarray*}%
Also, let $G_{n}$ be the semidirect product $G_{n}=V_{n}\left\langle
x_{n}\right\rangle $. Then $G_{n}$ is a $\mathfrak{T}_{n}$-group. Let $p$ be
a prime number, $W_{n}=pV_{n}$, $x^{\prime }=v_{n}x$ and $%
H_{n}=W_{n}\left\langle x_{n}^{\prime }\right\rangle $. Then, $\left[
G_{n}:H_{n}\right] =p^{n}$ and 
\begin{equation*}
f_{n}:pv_{i}\rightarrow v_{i}\text{ }\left( 1\leq i\leq n\right) ,\text{ }%
x_{n}^{\prime }\rightarrow x_{n}
\end{equation*}%
extends to a strongly simple epimorphism $f:H\rightarrow G$.

\subsection{Triples with degree a product of two primes}

Let $\left( G,H,f\right) $ be a simple triple of degree $m=pq$, a product of
two primes. In contrast to the prime degree, here we have a greater variety
of groups. However, we do not know of examples of simple $\left(
G,H,f\right) $ where $G$ is a non-abelian $\mathfrak{T}$-group, $m=pq$ and $%
p $, $q$ distinct primes. We show below that $G$ can be a mixed non-abelian
group.

\textbf{Examples 6. }(1)\textbf{\ }\textit{A mixed nilpotent group} $G$ 
\textit{of class} $2$, \textit{with} $[G:H]=p^{2}$.

Let $G=<a,u,v$ $|$ $u^{p}=v^{p}=\left[ u,v\right] =e,u^{a}=uv,v^{a}=v>$ of
type $\left( \mathbb{Z}_{p}\mathbb{\times Z}_{p}\right) \mathbb{Z}$. Let $%
H=<u,a^{p}>$ and $K=<v,a>$. Then $[G:H]=p^{2}$ and $H,K$ are abelian of type 
$\mathbb{Z}_{p}\mathbb{\times Z}$.

Define $f:H\rightarrow K$ by $f:u\rightarrow v,a^{p}\rightarrow a,$. Then, $%
f $ is a simple homomorphism; for%
\begin{equation*}
u^{a^{p}}=uv^{p}=u,f:u^{a^{p}}\left( =u\right) \rightarrow v^{a}\left(
=v\right) .
\end{equation*}

(2)\textbf{\ }\textit{A mixed nilpotent group} $G$ \textit{of class} $2$, 
\textit{with} $[G:H]=pq$ \textit{where} $p$ \textit{is a prime}, $q=1+tp$ 
\textit{and} $\left( G,H,f\right) $ \textit{simple}.

Let $D=\left\langle a,b|a^{p^{2}}=b^{p}=e,a^{b}=a^{1+p}\right\rangle $, a
group of order $p^{3}$. Let $G=D.\left\langle x\right\rangle $ where $x$ is
of infinite order and acts on $D$ as conjugation by $b$. Then,%
\begin{equation*}
Z\left( G\right) =\left\langle a^{p},b^{-1}x\right\rangle \text{.}
\end{equation*}%
We observe that $\left( b^{-1}x\right) ^{p}=$ $x^{p}$. Let $q=1+tp$ and let%
\begin{equation*}
H=\left\langle a^{p},b,x^{q}\right\rangle ,K=\left\langle b,x\right\rangle 
\text{. }
\end{equation*}%
Then, $H$ is abelian, has index $pq$ in $G$, has type $\mathbb{Z}_{p}\mathbb{%
\times Z}_{p}\mathbb{\times Z}$ and $K$ is abelian of type $\mathbb{Z}_{p}%
\mathbb{\times Z}$. Moreover, 
\begin{equation*}
H\cap Z\left( G\right) =\left\langle a^{p},b^{-1}x^{q}\right\rangle
\end{equation*}%
since $b^{-1}x.\left( x^{p}\right) ^{t}=b^{-1}x^{q}$. Now the map 
\begin{equation*}
f:a^{p}\rightarrow b,b\rightarrow e,x^{q}\rightarrow x
\end{equation*}%
extends to an epimorphism $f:H\rightarrow K$. To prove that $f$ is simple we
observe that 
\begin{eqnarray*}
f &:&H\cap Z\left( G\right) \rightarrow \left\langle b,x\right\rangle , \\
H\cap \left\langle b,x\right\rangle &=&\left\langle b,x^{q}\right\rangle , \\
f &:&\left\langle b,x^{q}\right\rangle \rightarrow \left\langle
x\right\rangle \text{.}
\end{eqnarray*}

\begin{theorem}
Let $G$ be a $\mathfrak{T}_{c}$-group, $\left( G,H,f\right) $ a simple$\ $of
degree $pq$ where $p,q$ are (not necessarily distinct) prime numbers and let 
$\frac{L}{Z(G)}=\overline{f}$-$core\left( \frac{HZ(G)}{Z(G)}\right) $. Then,%
\newline
(i) $L$ and $\frac{G}{L}$ are free abelian groups;\newline
(ii) $G=HH^{f}$;\newline
(iii) $Z_{c-1}(G)\leq L$;\newline
(iv) $Z(G)=\sqrt[G]{\gamma _{c}\left( G\right) }$.
\end{theorem}

\begin{proof}
Suppose $G$ is non-abelian. Since $Z(H)\not=Z(G)$, we may suppose $\left[
G:HZ(G)\right] =p,\left[ HZ(G):H\right] =$ $q$. Then, $\left[ Z(G):Z(H)%
\right] =q$ and $HZ(G)=HL$.

The subgroups $HZ(G)$, $\left( HZ(G)\right) ^{\prime }=H^{\prime }$ are
normal in $G$. Since $(\frac{G}{L},\frac{HL}{L},\overline{f})$ is a simple
triple of degree a prime, by Corollary 2, $\frac{G}{L}$ is either free
abelian or finite. We know that $L$ and $\sqrt[G]{L}$ are abelian.
Therefore, $\frac{G}{L}$ is free abelian.

As $Z(H)^{f}\not\leq Z(H)$ and $\left[ Z(G):Z(H)\right] =q$, we obtain $%
Z(G)=Z(H)Z(H)^{f}$. It follows from $\left( HZ(G)\right) H^{f}=HH^{f}$ that $%
HH^{f}$ is a subgroup of $G$ and $HH^{f}=$ $HZ(G)$ or $G$. Suppose the first
alternative holds. Then, $\left( HH^{f}\right) ^{\prime }=$ $H^{\prime }$, a
normal subgroup of $G$ and $(H^{f})^{\prime }=(H^{\prime })^{f}\leq $ $%
H^{\prime }$. Therefore, $H$ is abelian and hence, central; this leads to $%
\left[ G:Z(G)\right] =p$ which is absurd.

Consider the first index $j$ such that $Z_{j}(G)\not\leq L$ and define $%
K=Z_{j}(G)L$. We assert that $\gamma _{i}\left( K\right) \leq Z_{j-i+1}(G)$
for $i\geq 2$:%
\begin{equation*}
\gamma _{2}\left( K\right) =Z_{j}(G)^{\prime }\left[ L,Z_{j}(G)\right] \leq
Z_{j-1}(G)\leq L\text{;}
\end{equation*}%
if $\gamma _{i}\left( K\right) \leq Z_{j-i+1}(G)$ for some $i$ then 
\begin{eqnarray*}
\gamma _{i+1}\left( K\right) &\leq &\left[ Z_{j-i+1}(G),Z_{j}(G)L\right] \\
&\leq &\left[ Z_{j-i+1}(G),Z_{j}(G)\right] \left[ Z_{j-i+1}(G),L\right] \\
&\leq &Z_{j-i+2}(G)\text{.}
\end{eqnarray*}

Since $\left( L\cap H\right) ^{f}\leq L$, and $\left(
Z_{j}(G),Z_{j}(H),f\right) $ is a triple, we have the corresponding triple $(%
\frac{K}{L},\frac{LZ_{j}(H)}{L},\overline{f})$. Furthermore, as $\frac{G}{L}$
is free abelian and $(\frac{G}{L},\frac{LH}{L},\overline{f})$ is simple of
degree $p$, it follows that $\frac{K}{L}$ is of finite index in $\frac{G}{L}$
and therefore, $K$ and $G$ have the same nilpotency class; that is, $G$ has
class $j=c$. Hence, $G=Z_{c}(G)=K$ and $Z_{c-1}(G)\leq L$.

Since $\left( Z(G),Z(H),f\right) $ is a simple triple of prime degree and $%
\left( \gamma _{c}(G),\gamma _{c}(H),f\right) $ is a sub-triple, it follows
that $\left[ Z(G):\gamma _{c}(G)\right] $ is finite; hence $Z(G)=\sqrt[G]{%
\gamma _{c}(G)}$.
\end{proof}

\subsection{A sequence of simple triples of degree $4$}

The groups $R,S$ produced in Subsection 3.2 will be shown to be part of an
ascending sequence of simple triples $\left( G_{n},H_{n},f_{n}\right) $
where $\left[ G_{n}:H_{n}\right] =2^{2}$, $d\left( G_{n}\right) =2$, $%
s\left( G_{n}\right) =2$ and $c\left( G_{n}\right) =n$. This will prove that
the nilpotency class of groups in Theorem 12 cannot have a fixed upper limit.

Let $V_{n}$ be the free $\mathbb{Z}$ module $\mathbb{Z}^{n}$ and $\left\{
\varepsilon _{i}|1\leq i\leq n\right\} $ its canonical basis. Define
inductively $x_{n}\in GL\left( V_{n}\right) $: 
\begin{eqnarray*}
x_{2} &=&\left( 
\begin{array}{cc}
1 & 1 \\ 
0 & 1%
\end{array}%
\right) , \\
x_{n} &=&\left( 
\begin{array}{cc}
1 & \xi _{n} \\ 
0 & x_{n-1}%
\end{array}%
\right) ,\xi _{n}=\left( 2^{n-2},0,...,0\right)
\end{eqnarray*}%
for $n\geq 3$. Then $x_{n}$ acts nilpotently and uniserially on $V_{n}$. We
note that $x_{n}^{2}$ leaves invariant the submodule $W_{n}=\left\langle
\varepsilon _{1},..,\varepsilon _{n-2},\varepsilon _{n-1}+\varepsilon
_{n},2\varepsilon _{n}\right\rangle $, where clearly, $\left[ V_{n}:W_{n}%
\right] =2$. Define the semidirect product $G_{n}=V_{n}\left\langle
x_{n}\right\rangle $ and its subgroup $H_{n}=W_{n}\left\langle
x_{n}^{2}\right\rangle $. Then, $G_{n}$ is nilpotent of class $n$ and $%
Z\left( G_{n}\right) =\left\langle \varepsilon _{n}\right\rangle $. Also, $%
\left[ G_{n}:H_{n}\right] =2^{2}$ and $Z\left( H_{n}\right) =\left\langle
2\varepsilon _{n}\right\rangle $. Furthermore, $G_{2}\cong R,G_{3}\cong S$.

We construct inductively simple endomorphisms $f_{n}:W_{n}\rightarrow V_{n}$:

\begin{eqnarray*}
f_{2} &=&\left( 
\begin{array}{cc}
1 & \frac{1}{2} \\ 
0 & \frac{1}{2}%
\end{array}%
\right) , \\
f_{n} &=&\left( 
\begin{array}{cc}
2^{n-2} & \theta _{n} \\ 
0 & f_{n-1}%
\end{array}%
\right) ,\theta _{n}=\left( \theta _{n,1},...,\theta _{n,n-1}\right)
\end{eqnarray*}%
for $n\geq 3$.

We extend $f_{n}$, first by stipulating that $f_{n}:x_{n}^{2}\rightarrow
x_{n}$. Then, $f_{n}$ extends to a homomorphism $H_{n}\rightarrow $ $G_{n}$
if and only if $x_{n}^{2}f_{n}=f_{n}x_{n}$. It is easy to see that $f_{2}$
satisfies this equation. For $n\geq 3$, $x_{n}^{2}f_{n}=f_{n}x_{n}$ is
equivalent to 
\begin{equation*}
\theta _{n}\left( 1-x_{n-1}\right) =\xi _{n}\left( 2^{n-2}-\left(
1+x_{n-1}\right) f_{n-1}\right) \text{.}
\end{equation*}%
This last equation provides a unique solution $\theta _{n}$ where 
\begin{eqnarray*}
\theta _{n,1} &=&2^{2}\theta _{n-1,1}+2^{2n-6}\,\text{,} \\
\theta _{n,i} &=&2^{i+1}\left( \theta _{n-1,i}+2^{n-4}\theta
_{n-2,i-1}\right) \text{ }\left( 2\leq i\leq \left[ \frac{n}{2}\right]
\right)  \\
&=&0\text{ }\left( \lfloor \frac{n+1}{2}\rfloor \leq i\leq n-1\right) \text{.%
}
\end{eqnarray*}%
The referee provided the following explicit solution%
\begin{equation*}
\theta _{n,i}=\frac{n}{i}\QOVERD( ) {n-i-1}{i-1}2^{\left( n-\frac{i}{2}%
-3\right) \left( i+1\right) +1}\text{.}
\end{equation*}%
The first few vectors are%
\begin{eqnarray*}
\theta _{3} &=&\left( 3,0\right) ,\theta _{4}=\left( 2^{4},2^{2},0\right) ,
\\
\theta _{5} &=&\left( 2^{4}5,2^{4}5,0,0\right) ,\theta _{6}=\left(
2^{7}3,2^{7}3^{2},2^{8},0,0\right) \text{.}
\end{eqnarray*}%
Finally, the resulting $f_{n}$ is simple. For otherwise, if $K$ is a
non-trivial subgroup of $H_{n}$, which is normal in $G_{n}$ and is invariant
under $f_{n}$ then $\left\{ e\right\} \not=Z\left( H_{n}\right) \cap K$
would also be $f_{n}$-invariant.

\section{ Index Theorem and an Application}

Let $G$ be a $\mathfrak{T}_{c}$-group, $H$ a subgroup of finite index $m$.
and $f:H\rightarrow G$ simple. Then, $\left[ G:H^{f}\right] =m^{\prime }$, $%
\left[ \gamma _{c}\left( G\right) :\gamma _{c}\left( H\right) \right] =l,$ $%
\left[ \gamma _{c}\left( G\right) :\gamma _{c}\left( H^{f}\right) \right]
=l^{\prime }$ are finite. The following general result establishes a simple
arithmetic relation between $m,m^{\prime },l,l^{\prime }$.

\begin{theorem}
Let $G$ be finitely generated nilpotent group, $H$ a subgroup of $G$ of
finite index $\left[ G:H\right] =m$ and $f:H\rightarrow G$ a monomorphism,
also let $\left[ G:H^{f}\right] =m^{\prime }$. Furthermore, let $U$ be a
subgroup of $H$ and let $V=\left\langle U,U^{f}\right\rangle $. Suppose $%
\left[ V:U\right] =l,\left[ V:U^{f}\right] =l^{\prime }$ are finite. Then
there exist integers $m_{1}|m,m_{1}^{\prime }|m^{\prime }$ such that $%
lm_{1}^{\prime }=l^{\prime }m_{1}$.
\end{theorem}

Before giving the proof, we note that the theorem is clearly true for finite
groups, since $\left\vert H\right\vert =\left\vert H^{f}\right\vert
,\left\vert U\right\vert =\left\vert U^{f}\right\vert $. However, the
following example shows that the theorem is not valid for the class of $2$%
-generated metabelian groups.

\textbf{Example 7. }Let $p$ be a prime number and $P$ be the subgroup of the
additive rationals \ generated by $\left\{ p^{i}|i\in \mathbb{Z}\right\} $.
Then $P$ admits the automorphism $f:p^{i}\rightarrow p^{i-1}$. Define $G$ to
be the extension of $P$ by $\left\langle f\right\rangle $; then $G$ \ is
generated by $\left\{ 1,f\right\} $. Now let $H=G$ and $U=\left\langle
p\right\rangle \leq P$. Then, $U^{f}=\left\langle 1\right\rangle =V$.
Moreover, the indices are $\left[ G:H\right] =1=\left[ G:H^{f}\right] ,\left[
V:U\right] =p,\left[ V:U^{f}\right] =1$.

\begin{proof}
I. Suppose $G$ is a free additive abelian group.

Then, 
\begin{equation*}
\sqrt[G]{V}=\sqrt[G]{U}=\sqrt[G]{U^{f}}
\end{equation*}%
and $\left[ \sqrt[G]{U}:U\right] =t,\left[ \sqrt[G]{U^{f}}:U^{f}\right]
=t^{\prime }$ are finite. Indeed, $t=t^{\prime }$, as $f$ induces an
isomorphism between the quotient groups $\frac{\sqrt[H]{U}}{U}$,$\frac{\sqrt[%
H^{f}]{U^{f}}}{U^{f}}$.

Define the indices 
\begin{equation*}
\left[ H+\sqrt[G]{V}:H\right] =m_{1},\left[ H^{f}+\sqrt[G]{V}:H^{f}\right]
=m_{1}^{\prime }
\end{equation*}%
where $m_{1}|m$ and $m_{1}^{\prime }|m^{\prime }$. Since 
\begin{equation*}
H\cap \sqrt[G]{V}=\sqrt[H]{U},H^{f}\cap \sqrt[G]{V}=\sqrt[H^{f}]{U^{f}}
\end{equation*}%
we conclude%
\begin{equation*}
\left[ \sqrt[G]{V}:\sqrt[H]{U}\right] =m_{1},\left[ \sqrt[G]{V}:\sqrt[H^{f}]{%
U^{f}}\right] =m_{1}^{\prime }\text{.}
\end{equation*}

Now, we calculate the index%
\begin{equation*}
\left[ \sqrt[G]{V}:U\cap U^{f}\right] =\left[ \sqrt[G]{V}:\sqrt[H]{U}\right] %
\left[ \sqrt[H]{U}:U\right] \left[ U:U\cap U^{f}\right] =
\end{equation*}%
\begin{equation*}
\left[ \sqrt[G]{V}:\sqrt[H^{f}]{U^{f}}\right] \left[ \sqrt[H^{f}]{U^{f}}%
:U^{f}\right] \left[ U^{f}:U\cap U^{f}\right] \text{.}
\end{equation*}%
Thus, $m_{1}tl^{\prime }=m_{1}^{\prime }tl$ and we reach $m_{1}l^{\prime
}=m_{1}^{\prime }l$.

We have $Tor\left( H\right) ^{f}=Tor\left( H^{f}\right) \leq Tor(G)$. The
map $\overline{f}:\frac{Tor(G)H}{Tor(G)}\rightarrow \frac{G}{Tor(G)}$ where $%
Tor(G)h\rightarrow Tor(G)h^{f}$ is a well-defined monomorphism. Since $\frac{%
G}{Tor(G)}$ is torsion-free nilpotent, it follows that $\left[ \frac{G}{%
Tor(G)}:\frac{Tor(G)H^{f}}{Tor(G)}\right] $ is finite and therefore $\left[
G:H^{f}\right] =m^{\prime }$ is finite.

II. Suppose $Tor\left( G\right) =\left\{ e\right\} $ We proceed by induction
on the nilpotency class of $G$. The case where $G$ is free-abelian was done
in part I.

We introduce the following notation 
\begin{eqnarray*}
Z_{H} &=&H\cap Z\left( G\right) ,Z_{H^{f}}=H^{f}\cap Z\left( G\right) \\
\left[ G:HZ\right] &=&m_{2},\left[ HZ:H\right] =m_{1},\left[ G:H^{f}Z\right]
=m_{2}^{\prime },\left[ H^{f}Z:H^{f}\right] =m_{1}^{\prime }\text{.}
\end{eqnarray*}

Then, 
\begin{eqnarray*}
Z\left( H\right) ^{f} &=&\left( H\cap Z\left( G\right) \right) ^{f}=\left(
Z_{H}\right) ^{f}, \\
\left[ Z:Z_{H}\right] &=&m_{1},\left[ Z:Z_{H^{f}}\right] =m_{1}^{\prime }%
\text{, }m=m_{2}m_{1},m^{\prime }=m_{2}^{\prime }m_{1}^{\prime }\text{.}
\end{eqnarray*}

Similarly, with respect to $V$, denote 
\begin{equation*}
Z_{V}=V\cap Z,Z_{U}=U\cap Z,Z_{U^{f}}=U^{f}\cap Z\text{,}
\end{equation*}%
and 
\begin{equation*}
\left[ V:UZ_{V}\right] =l_{2},\left[ UZ_{V}:U\right] =l_{1}\text{, }\left[
V:U^{f}Z_{V}\right] =l_{2}^{\prime },\left[ U^{f}Z_{V}:U^{f}\right]
=l_{1}^{\prime }\text{.}
\end{equation*}%
Then 
\begin{equation*}
\left[ Z_{V},Z_{U}\right] =l_{1},\left[ Z_{V},Z_{U^{f}}\right]
=l_{1}^{\prime }\text{, }l=l_{2}l_{1},l^{\prime }=l_{2}^{\prime
}l_{1}^{\prime }\text{.}
\end{equation*}

We claim that $\left( Z_{U}\right) ^{f}=Z_{U^{f}}$. This follows from 
\begin{equation*}
\left( Z_{U}\right) ^{f}\leq Z\left( H^{f}\right) \leq Z\left( G\right)
,\left( Z_{U}\right) ^{f}\leq U^{f}\cap Z\left( G\right) =Z_{U^{f}}
\end{equation*}%
and from 
\begin{eqnarray*}
Z_{U^{f}} &\leq &Z\left( H^{f}\right) =Z\left( H\right) ^{f}=\left( Z\left(
G\right) \cap H\right) ^{f}, \\
Z_{U^{f}} &\leq &\left( Z\left( G\right) \cap H\right) ^{f}\cap U^{f}=\left(
Z\left( G\right) \cap U\right) ^{f}=\left( Z_{U}\right) ^{f}\text{.}
\end{eqnarray*}

From the configuration 
\begin{equation*}
Z\geq Z_{H},Z_{H^{f}}\text{, }Z_{V}\geq Z_{U},Z_{U^{f}}
\end{equation*}%
we obtain that there exist $m_{11}|m_{1},m_{11}^{\prime }|m_{1}^{\prime }$
such that $l_{1}m_{11}^{\prime }=l_{1}^{\prime }m_{11}$.

Next, we apply induction to the nilpotency class of $G$.

We have $Z(H)\leq Z\left( G\right) $ and $Z(H)^{f}=Z\left( H^{f}\right) \leq
Z\left( G\right) $. Moreover, $\overline{f}:\frac{Z\left( G\right) H}{%
Z\left( G\right) }\rightarrow \frac{G}{Z\left( G\right) }$ defined by $%
Z\left( G\right) h\rightarrow Z\left( G\right) h^{f}$ is a monomorphism and $%
(\frac{Z\left( G\right) U}{Z\left( G\right) })^{\overline{f}}=\frac{Z\left(
G\right) U^{f}}{Z\left( G\right) }$. By applying induction to the class of $%
\frac{G}{Z\left( G\right) }$, we obtain that there exist $%
m_{21}|m_{2},m_{21}^{\prime }|m_{2}^{\prime }$ such that $%
l_{2}m_{21}^{\prime }=l_{2}^{\prime }m_{21}$. Hence, putting together the
two equations $l_{1}m_{11}^{\prime }=l_{1}^{\prime
}m_{11},l_{2}m_{21}^{\prime }=l_{2}^{\prime }m_{21}$, we obtain

\begin{equation*}
m_{1}^{\prime }=m_{11}^{\prime }m_{21}^{\prime }|m^{\prime
},m_{1}=m_{11}m_{21}|m
\end{equation*}%
\begin{eqnarray*}
l_{1}m_{11}^{\prime }l_{2}m_{21}^{\prime } &=&l_{1}^{\prime
}m_{11}l_{2}^{\prime }m_{21}, \\
lm_{1}^{\prime } &=&l^{\prime }m_{1}\text{.}
\end{eqnarray*}

III. Now we argue the general case where $T=Tor(G)$ is not necessarily
trivial.

Similar to the work done in part II, we define $T_{H}=T\cap H$ and likewise
we define

$T_{H^{f}}$, $T_{V},T_{U},T_{U^{f}}$, all finite groups. We note also that $%
\left( T_{U}\right) ^{f}=T_{U^{f}}$ and therefore $T_{U},T_{U^{f}}$ have
equal orders. Then it follows that 
\begin{equation*}
\left[ V:UT_{V}\right] =\left[ T_{V}:T_{U}\right] =\left[ T_{V}:T_{U^{f}}%
\right] =\left[ V:U^{f}T_{V}\right] \text{.}
\end{equation*}%
Finally, the argument continues as in part (II) with $T$ substituting $Z$.
\end{proof}

A special case of the above result is

\begin{corollary}
Maintain the hypotheses of the theorem. (i) If $U^{f}\leq U$ then $l=\left[
U:U^{f}\right] $ is finite and $l|m$. (ii) If $U\leq U^{f}$ then $l^{\prime
}=\left[ U^{f}:U\right] $ is finite and $l^{\prime }|m^{\prime }$.
\end{corollary}

\begin{proof}
Suppose $U^{f}\leq U$. Then as $U^{f}$ is isomorphic to $U$, we conclude $%
\left[ U:U^{f}\right] $ is finite. The remaining assertions are direct.
\end{proof}

We use the above divisibility criterion to prove

\begin{theorem}
Let $G$ be a $\mathfrak{T}_{c}$-group, $H$ a subgroup of finite index $m$ in 
$G$ and $f:H\rightarrow G$ an epimorphism. Let $\left[ G:HZ_{c-1}\left(
G\right) \right] =$ $k$, $\left[ Z\left( G\right) :Z\left( H\right) \right]
=q$. Then, $\left[ \gamma _{c}\left( G\right) :\gamma _{c}\left( H\right) %
\right] $ is a $k$-number which divides $q$. If $f$ is simple and $m$ a
square-free integer then $G$ is abelian.
\end{theorem}

\begin{proof}
We may assume $c\geq 2$. By Subsection 3.1, item I.2, if $G=HZ_{c-1}\left(
G\right) $ then $\gamma _{c}\left( G\right) =$ $\gamma _{c}\left( H\right) $%
. So, suppose $G\not=HZ_{c-1}\left( G\right) $ and consider the free abelian
group $\overline{G}=\frac{G}{Z_{c-1}\left( G\right) }$. There exist $%
a_{1},a_{2},...,a_{s}\in G$ such that%
\begin{equation*}
Z_{c-1}\left( G\right) a_{1},Z_{c-1}\left( G\right) a_{2},...,Z_{c-1}\left(
G\right) a_{s}
\end{equation*}%
freely generate $\overline{G}$ and integers $k_{1}|k_{2}|...|k_{r}$, $%
k=k_{1}k_{2}...k_{r}$ such that%
\begin{equation*}
Z_{c-1}\left( G\right) a_{1}^{k_{1}},Z_{c-1}\left( G\right)
a_{2}^{k_{2}},...Z_{c-1}\left( G\right) a_{r}^{k_{r}},Z_{c-1}\left( G\right)
a_{r+1}..,Z_{c-1}\left( G\right) a_{s}
\end{equation*}%
freely generate $\overline{H}=\frac{HZ_{c-1}\left( G\right) }{Z_{c-1}\left(
G\right) }$. Thus, there exist $c_{1},c_{2},...,c_{s}\in Z_{c-1}\left(
G\right) $ such that%
\begin{equation*}
H=\left\langle
c_{1}a_{1}^{k_{1}},c_{2}a_{2}^{k_{2}},...,c_{r}a_{r}^{k_{r}},c_{r+1}a_{r+1}..,c_{s}a_{s}\right\rangle Z_{c-1}\left( H\right) 
\text{.}
\end{equation*}

Now, $\gamma _{c}\left( G\right) \leq $ $Z\left( G\right) $ and is generated
by simple commutators $\left[ a_{i_{1}},a_{i_{2}},...,a_{i_{c}}\right] $ of
weight $c$ where the indices $i_{j}$ are from $\left\{ 1,2,...,s\right\} $.
Whereas, $\gamma _{c}\left( H\right) $ is generated by $\left[
a_{i_{1}},a_{i_{2}},...,a_{i_{c}}\right] ^{\lambda \left(
i_{1},...,i_{c}\right) }$ where $\lambda \left( i_{1},...,i_{c}\right)
=k_{1}^{u_{1}}...k_{t}^{u_{t}}\not=1$ and $u_{z}$ is the number of $%
i_{j}=z\in \left\{ 1,...,r\right\} $. Therefore, $\left\vert \frac{\gamma
_{c}\left( G\right) }{\gamma _{c}\left( H\right) }\right\vert $ is a $k$%
-number.

As $f$ induces epimorphisms $\gamma _{c}\left( H\right) \rightarrow \gamma
_{c}\left( G\right) $, $Z\left( H\right) \rightarrow Z\left( G\right) $, we
apply Corollary 3 to obtain $\left\vert \frac{\gamma _{c}\left( G\right) }{%
\gamma _{c}\left( H\right) }\right\vert $ divides $\left\vert \frac{Z\left(
G\right) }{Z\left( H\right) }\right\vert =q$.

Now suppose $f$ is simple, $m$ square-free and $G$ non-abelian. Since $f$ is
simple, we have $Z_{c-1}\left( G\right) \not\leq H$. Let%
\begin{equation*}
\left[ G,HZ_{c-1}\left( G\right) \right] =m_{1},\left[ HZ_{c-1}\left(
G\right) ,H\right] =m_{2}\text{;}
\end{equation*}%
then, $\gcd \left( m_{1},m_{2}\right) =1$. Since $\left[ Z\left( G\right)
:Z\left( H\right) \right] $ divides $m_{2}$, we conclude that $\gamma
_{c}\left( G\right) =\gamma _{c}\left( H\right) $ and therefore, $\gamma
_{c}\left( G\right) =\left\{ e\right\} $; a contradiction.
\end{proof}

\end{document}